\documentclass[twocolumn]{jdjnmacs}
\usepackage[latin1]{inputenc}
\usepackage[T1]{fontenc}
\usepackage[francais]{babel}
\usepackage{graphicx}
\usepackage{psfrag}
\usepackage{amsfonts}
\usepackage{amssymb,amsbsy,amsmath,amsfonts,amssymb,amscd}

\newtheorem{remarque}{\it Remarque\/}

\begin{document}

\title{\vspace{0.5cm}Rien de plus utile qu'une bonne th\'{e}orie: \\ la commande sans mod\`{e}le
 \\ ~ \\ {\em Nothing is as Practical as a Good Theory$:$
\\ Model-Free Control}}
\author{\vskip 1em Michel \textsc{FLIESS}\textsuperscript{1}, C\'{e}dric \textsc{JOIN}%
\textsuperscript{2, 4}, Samer \textsc{RIACHY}\textsuperscript{3,  4} \\
\vskip1em \textsuperscript{1} LIX (CNRS, UMR 7161), \'Ecole
polytechnique, 91228 Palaiseau, France\\
\vskip 1em \textsuperscript{2} CRAN (CNRS, UMR 7039),
Nancy-Universit\'{e}, BP 239, 54506 Vand\oe{}uvre-l\`es-Nancy, France\\
\vskip 1em \textsuperscript{3} ECS (EA 3649), ENSEA, 6 avenue du
Ponceau, 95014 Cergy-Pontoise, France\\ \vskip 1em
\textsuperscript{4}\'Equipe Non-A, INRIA  Lille -- Nord-Europe
\vskip 1em \texttt{Michel.Fliess@polytechnique.edu,
Cedric.Join@cran.uhp-nancy.fr, \\
Samer.Riachy@ensea.fr} }

\maketitle

\begin{abstract}
L'exp\'{e}rience acquise avec de nombreuses applications r\'{e}ussies permet
de revisiter la commande sans mod\`{e}le. On remplace la d\'{e}rivation
num\'{e}rique de signaux bruit\'{e}s par une identification param\'{e}trique en
temps r\'{e}el, beaucoup plus simple. On explique la curieuse
universalit\'{e} des PID usuels et pourquoi un mod\`{e}le ultra-local
d'ordre 1 suffit presque toujours. On montre que, m\^{e}me avec un
mod\`{e}le restreint, l'usage d'un correcteur PI intelligent reste
avantageux. Deux exemples, un pendule invers\'{e} et un \'{e}changeur de
vapeur, semblent confirmer l'avantage de la commande sans-mod\`{e}le par
rapport \`{a} d'autres approches, comme la commande par r\'{e}gimes
glissants et l'utilisation de techniques de dimension infinie, qu'il
s'agisse de retards ou d'\'{e}quations aux d\'{e}riv\'{e}es partielles.
Plusieurs simulations num\'{e}riques illustrent notre
propos. 
\newline \textit{Abstract}--- The experience gained with numerous
successful applications permits to revisit some points of model-free
control. The numerical differentiation of noisy signals may be
replaced by a real time parameter identification which is much
simpler. The strange ubiquity of classic PIDs is explained as well
as the almost universal utilization of ultra-local models of order
1. We show that even with a partially known model the utilization of
an intelligent PI controller remains profitable. Two examples, an
inverted pendulum and a heat exchanger, seem to confirm the
superiority of model-free control with respect to sliding modes
control and techniques stemming from infinite-dimensional systems,
such as delays or partial differential equations. Our paper is
illustrated by several computer simulations. \\ ~ \\ ~ \\
\end{abstract}

\begin{keywords}
Commande sans mod\`{e}le, PI intelligent, PID, identification, pendule
invers\'{e}, \'{e}changeur de chaleur, r\'{e}gimes glissants, retards, \'{e}quations
aux d\'{e}riv\'{e}es partielles. \\
\textit{Keywords}--- Model-free
control, intelligent PI, PID, identification, inverted pendulum,
heat exchanger, sliding modes, delays, partial differential
equations.
\end{keywords}
\newpage

\section{Introduction}
\og Rien de plus utile qu'une bonne th\'{e}orie\footnote{Traduite de
l'anglais: \og \textit{nothing is as practical as a good theory}
\fg, cette phrase est due au psychologue Kurt Lewin. Les auteurs
remercient chaleureusement Roger Tauleigne (ECS, ENSEA,
Cergy-Pontoise) de la leur avoir sugg\'{e}r\'{e}e comme titre.} \fg, voil\`{a}
qui r\'{e}sume les ambitions de la \emph{commande sans mod\`{e}le}, due \`{a}
deux des auteurs (\cite{esta,malo}). Cette vision, radicalement
nouvelle de l'automatique\footnote{Voir, cependant,
\cite{chang-jung,han}.}, combine, rappelons-le, les deux avantages
suivants:
\begin{enumerate}
\item Inanit\'{e}, comme son nom l'indique, d'une
mod\'{e}lisation math\'{e}matique pr\'{e}cise, qu'elle soit bas\'{e}e sur des lois
physiques ou des proc\'{e}dures d'identification, comme l'exige trop
souvent la th\'{e}orie du contr\^{o}le \og moderne \fg. Frottements,
hyst\'{e}r\'{e}sis, effets thermiques, vieillissement, dispersion des
caract\'{e}ristiques due \`{a} la fabrication en s\'{e}rie, $\dots$, sont, on ne
le sait que trop, rebelles \`{a} toute description simple et fiable par
\'{e}quations diff\'{e}rentielles.
\item Facilit\'{e} du r\'{e}glage du correcteur correspondant, dit \emph{PI intelligent}, ou
\emph{iPI}. D'o\`{u} une rupture nette avec les PID classiques dont le
r\'{e}glage, malais\'{e} et p\'{e}nible (voir, par exemple, \cite{astrom,od}),
rel\`{e}ve trop souvent de \og recettes de cuisine \fg.
\end{enumerate}
De nombreuses applications (voir
\cite{agadir,cifa,buda,nolcos,brest,psa,edf,michel,vil1,vil2,vil3,vil4})
ont d\'{e}j\`{a} \'{e}t\'{e} r\'{e}ussies. Cet article tire parti de l'exp\'{e}rience
acquise pour am\'{e}liorer la pr\'{e}sentation.

Voici les modifications essentielles:
\begin{enumerate}
\item On remplace, comme en \cite{milan1}, la d\'{e}rivation
de signaux bruit\'{e}s (\cite{easy,mboup}) par l'identification
param\'{e}trique lin\'{e}aire, d\'{e}velopp\'{e}e en \cite{sira1,sira2},
conceptuellement plus simple, qui a suscit\'{e} par ailleurs quelques
applications (\cite{abou,becedas,spain,trapero}).
\item On explique
\begin{itemize}
\item pourquoi un mod\`{e}le \emph{ultra-local} d'ordre $1$ suffit presque toujous
en pratique;
\item l'\'{e}trange ubiquit\'{e} des PID classiques, malgr\'{e} leurs
r\'{e}glages fastidieux (voir \cite{pid}).
\end{itemize}
\item On montre, suivant \cite{milan1} et contrairement \`{a}
\cite{esta,malo}, qu'avec un mod\`{e}le restreint, c'est-\`{a}-dire
partiellement connu, on peut conserver un correcteur iPI.
\end{enumerate}
Apr\`{e}s des rappels au {\S} \ref{samo}, on traite ces questions aux {\S}
\ref{isch}, \ref{ord1}, \ref{connect} et \ref{nlspring}. Les
illustrations num\'{e}riques du {\S} \ref{illu}, c'est-\`{a}-dire un pendule
invers\'{e} et un \'{e}changeur de chaleur, mettent en \'{e}vidence la probable
sup\'{e}riorit\'{e} de notre approche par rapport aux modes glissants (voir,
aussi, \cite{milan1,milan2}) et questionne l'utilit\'{e} de techniques
issues de la dimension infinie, comme les retards (voir, aussi,
\cite{esta,malo}) et les \'{e}quations aux d\'{e}riv\'{e}es partielles (voir,
aussi, \cite{edf} et \cite{agadir}). La conclusion du {\S}
\ref{conclusion} questionne le futur de la recherche universitaire
en automatique et, plus g\'{e}n\'{e}ralement, la place pr\'{e}pond\'{e}rante de la
mod\'{e}lisation math\'{e}matique en sciences appliqu\'{e}es.

\section{Rappels sur le sans-mod\`{e}le}\label{samo}
\subsection{Le mod\`{e}le ultra-local}
On se restreint, pour simplifier les notations, \`{a} un syst\`{e}me
monovariable, d'entr\'{e}e $u$ et de sortie $y$. Dans l'ignorance d'un
mod\`{e}le math\'{e}matique global, on introduit le mod\`{e}le \og
ph\'{e}nom\'{e}nologique \fg, dit \emph{ultra-local}, valable sur un court
laps de temps,
\begin{equation}\label{F}
\boxed{y^{(\nu)} = F + \alpha u}
\end{equation}
o\`{u}
\begin{itemize}
\item l'ordre de d\'{e}rivation $\nu$, en g\'{e}n\'{e}ral $1$, choisi par
l'op\'{e}rateur, est \'{e}tranger \`{a} l'ordre de d\'{e}rivation maximum de $y$,
inconnu, dans le syst\`{e}me;

\item le param\`{e}tre constant $\alpha$, fix\'{e} par l'op\'{e}rateur afin que
les valeurs num\'{e}riques de $\alpha u$ et $y^{(\nu)}$ aient m\^{e}me ordre
de grandeur, n'a pas \textit{a priori} de valeur pr\'{e}cise;

\item $F$, qui contient toutes les informations \og structurelles
\fg, d\'{e}pend de toutes les autres variables du syst\`{e}me, y compris des
perturbations, et de leurs d\'{e}riv\'{e}es.
\end{itemize}
L'estimation en temps r\'{e}el de la valeur num\'{e}rique de $F$, trait\'{e}e au
{\S} \ref{oe}, permet de r\'{e}actualiser \eqref{F} \`{a} chaque instant.

\subsection{Correcteurs PI intelligents}
On obtient le comportement d\'{e}sir\'{e} avec, si $\nu = 1$ en \eqref{F},
le correcteur \emph{proportionnel-int\'{e}gral intelligent}, ou
\emph{iPI},
\begin{equation}\label{ipi}
\boxed{u =  - \frac{F - \dot{y}^\star + K_P e + K_I \int e}{\alpha}}
\end{equation}
o\`{u}
\begin{itemize}
\item $y^\star$ est la trajectoire de r\'{e}f\'{e}rence de la sortie,
\item $e = y - y^\star$ est l'erreur de poursuite,
\item $K_P$, $K_I$ sont les gains usuels.
\end{itemize}
$K_I = 0$ conduit au correcteur \emph{proportionnel intelligent}, ou
\emph{iP},
\begin{equation}\label{iP}
\boxed{u = - \frac{F - \dot{y}^\ast + K_P e }{\alpha}}
\end{equation}

\begin{remarque}
Si $\nu = 1$ en \eqref{F}, on se ram\`{e}ne avec \eqref{ipi} ou
\eqref{iP} \`{a} la stabilisation d'un int\'{e}grateur pur. D'o\`{u} le r\'{e}glage
facile des gains.
\end{remarque}

\section{Mise en {\oe}uvre}\label{isch}
\subsection{Estimation de $F$}\label{oe}
R\'{e}\'{e}crivons \eqref{ipi} sous la forme
\begin{equation}\label{est1}
F = - \alpha u + \dot{y}^\star - K_P e - K_I \int e
\end{equation}
On att\'{e}nue les bruits entachant les mesures en
int\'{e}grant\footnote{Les bruits, consid\'{e}r\'{e}s comme des fluctuations
rapides, sont att\'{e}nu\'{e}s par des filtres passe-bas, dont l'int\'{e}grale
est un exemple simple. Voir \cite{bruit} pour une explication
math\'{e}matique.} les deux membres de \eqref{est1} sur un court laps de
temps $\delta$. Il vient:
\begin{equation*}
F_{\text{\tiny approx}} = \frac{1}{\delta}\int_{T-\delta}^{T}\left (-\alpha u + \dot{y}^\star
- K_P e -K_I \int e\right ) d\tau \label{estimf}
\end{equation*}
o\`{u} $F_{\text{\tiny approx}}$ est une approximation, constante par
morceaux, de $F$. Cet estimateur s'implante facilement sous la forme
d'un filtre lin\'{e}aire discret.

\begin{remarque}
Il y a quelques situations o\`{u} il peut \^{e}tre int\'{e}ressant d'estimer
$\alpha$ en \eqref{F}. Voir, \`{a} ce sujet, \cite{sympl,milan2}.
\end{remarque}

\subsection{Trajectoire de r\'{e}f\'{e}rence pour la sortie}
On suppose le syst\`{e}me d'entr\'{e}e $u$ et de sortie $y$ \`{a} d\'{e}phasage
minimal\footnote{L'ignorance des \'{e}quations gouvernant le syst\`{e}me
interdit toute v\'{e}rification math\'{e}matique de la nature du d\'{e}phasage.
L'appr\'{e}ciation de cette nature repose donc sur une bonne
connaissance exp\'{e}rimentale du comportement. Voir \cite{pvtol} pour
une premi\`{e}re approche du d\'{e}phasage non minimal dans le cadre du
sans-mod\`{e}le.}. On peut alors choisir une trajectoire de r\'{e}f\'{e}rence
pour $y$ satisfaisant les n\'{e}cessit\'{e}s du syst\`{e}me.

\section{Pourquoi l'ordre $1$ suffit-il?}\label{ord1}
\begin{figure}[htp]
\centering
\includegraphics[width=9.05cm]{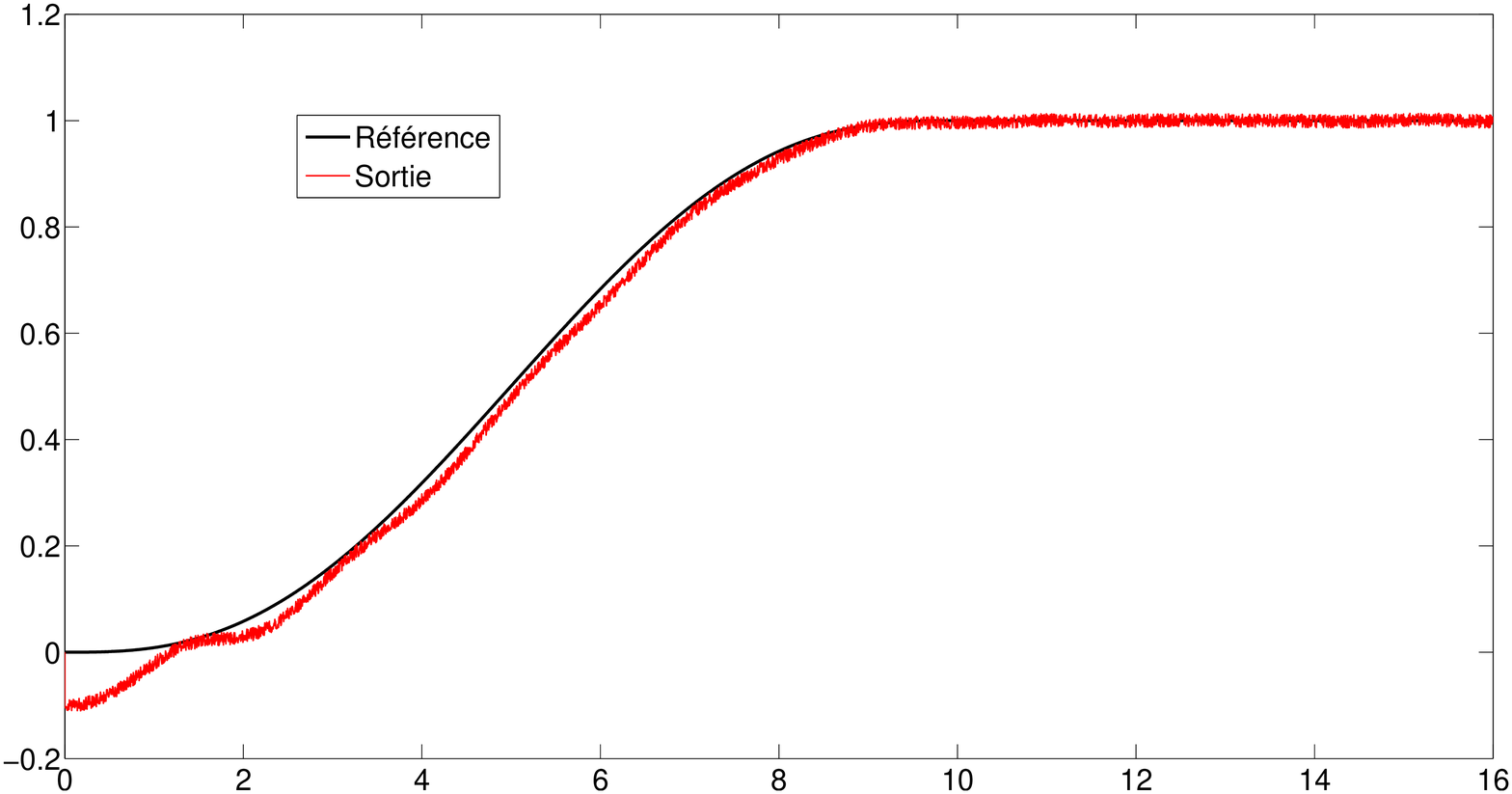}
\caption{Syst\`{e}me du 2$^e$ ordre avec frottement et r\'{e}gulateur
iPI.}\label{fig01}
\end{figure}

\begin{figure}[htp]
\centering
\includegraphics[width=9.05cm]{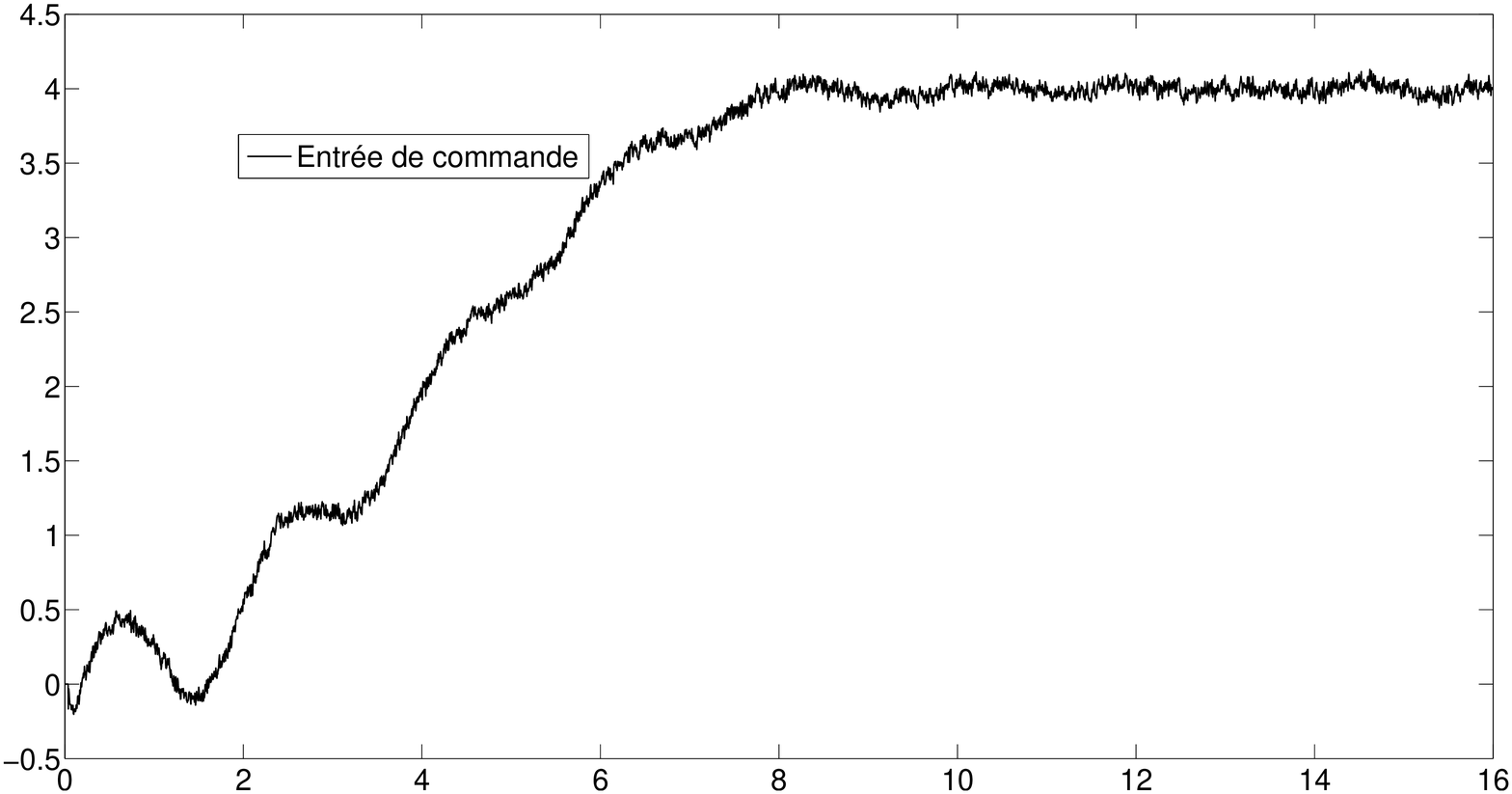}
\caption{Commande du r\'{e}gulateur iPI.}\label{fig02}
\end{figure}

\begin{figure}[htp]
\centering
\includegraphics[width=9.05cm]{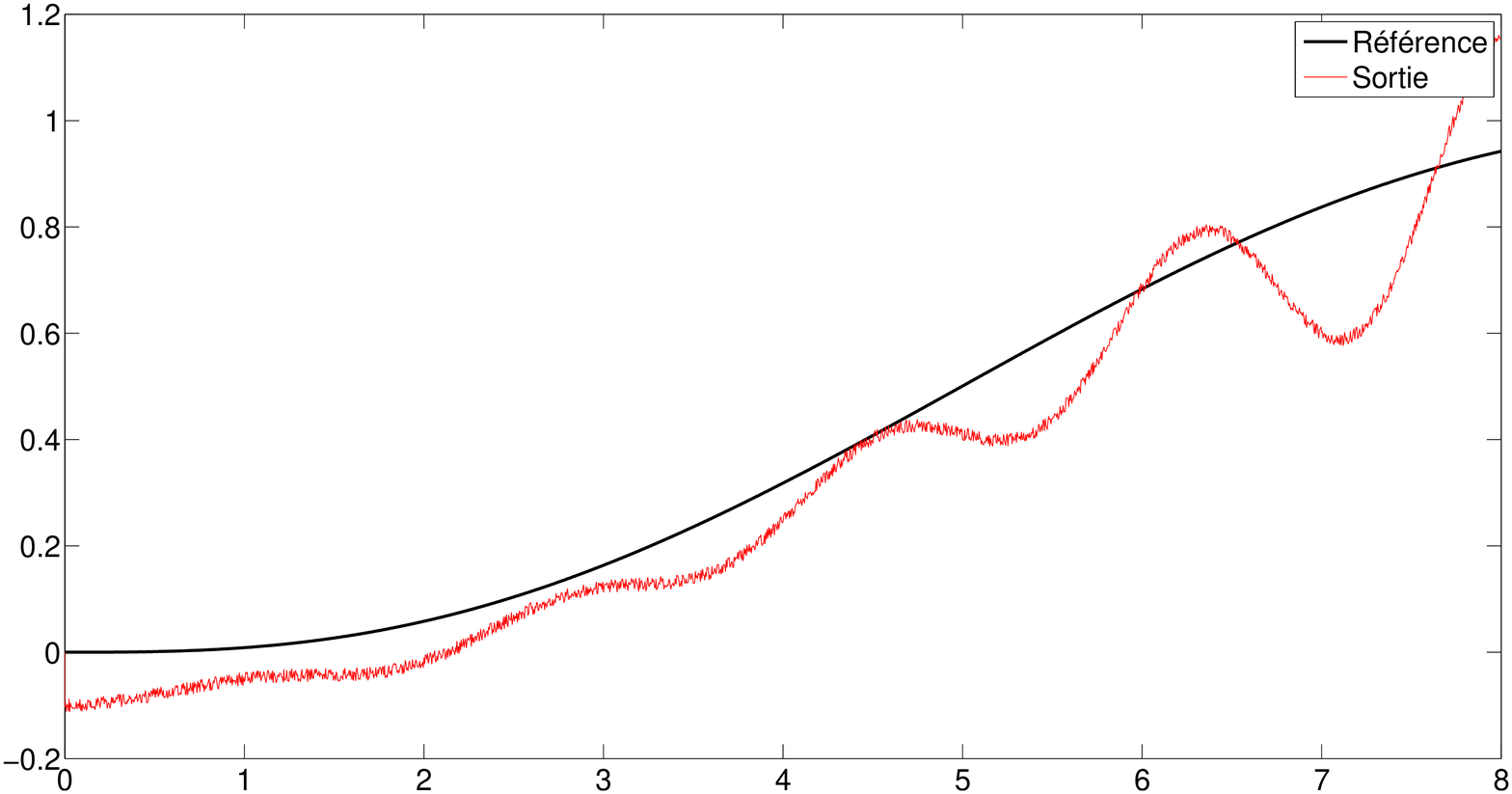}
\caption{Oscillateur harmonique avec r\'{e}gulateur iPI.}\label{fig03}
\end{figure}

Prendre $\nu = 1$ en \eqref{F} suffit jusqu'\`{a} pr\'{e}sent dans toutes
les applications\footnote{Sauf en \cite{buda}, o\`{u} $\nu = 0$.}.
Pourquoi cette heureuse propri\'{e}t\'{e} qui simplifie passablement
l'implantation num\'{e}rique du sans-mod\`{e}le? On l'explique par les
frottements. Il leur correspond la pr\'{e}sence de la d\'{e}riv\'{e}e premi\`{e}re
$\dot{y}$ dans l'\'{e}quation, inconnue, du syst\`{e}me, qui \'{e}vite
l'apparition d'une boucle alg\'{e}brique avec $\nu = 1$. Illustrons
cette th\`{e}se avec le syst\`{e}me lin\'{e}aire du second ordre, \`{a} coefficients
constants,
\begin{equation}
\ddot{x}+c\dot{x}+4 x=u \label{oscillateur}
\end{equation}
o\`{u} $c\dot{x}$ repr\'{e}sente des frottements. Les figures \ref{fig01},
\ref{fig02} fournissent des simulations satisfaisantes avec un
r\'egulateur iPI; le choix des valeurs des param\`{e}tres est $c=3$,
$\alpha = 1$, $K_P=16$, $K_I=25$. Une d\'{e}gradation consid\'{e}rable des
performances, comme on le voit sur la figure \ref{fig03}, se produit
avec un oscillateur harmonique, c'est-\`{a}-dire avec $c = 0$ en
\eqref{oscillateur}.


\section{Efficacit\'{e} des PI classiques}\label{connect}
\subsection{Discr\'{e}tisation}\label{1}
Associons au correcteur PI classique
\begin{equation}\label{cpi}
  u (t) = k_p e(t) + k_i \int e(\tau) d\tau
\end{equation}
sa forme \og vitesse \fg:
$$
\dot{u} (t) = k_p \dot{e}(t) + k_i e(t)
$$
Une discr\'{e}tisation \'{e}l\'{e}mentaire fournit:
\begin{equation}
\label{PId} \dfrac{u(t) - u(t - h)}{h} =  k_p \left(\dfrac{e(t) -
e(t - h)}{h}\right) + k_i  e(t)
\end{equation}
o\`{u} $h > 0$ est \og petit \fg.

Si $\nu = 1$ en \eqref{F}, rempla\c{c}ons $F$ par ${\dot y}(t) - \alpha
u (t-h)$ pour les besoins de l'implantation num\'{e}rique. Alors,
\eqref{iP} devient
\begin{equation}
\label{iPd} u (t) = u (t - h) - \frac{e(t) - e(t-h)}{h\alpha} +
\dfrac{K_P}{\alpha}\, e(t)
\end{equation}
\eqref{PId} et \eqref{iPd} deviennent identiques si, et seulement
si,
\begin{align}
\label{eqPI_i-P_corresp} k_p &= - \dfrac{1}{\alpha h}, \quad k_i =
\dfrac{K_P}{\alpha h}
\end{align}

\begin{remarque}
Cette \'{e}quivalence entre PI et iP n'est plus valable en temps
continu, comme on le voit si $h \rightarrow 0$.
\end{remarque}

\subsection{Explication}
Les calculs pr\'{e}c\'{e}dents, qui s'\'{e}tendent ais\'{e}ment aux PID (voir
\cite{pid}), expliquent pourquoi, gr\^{a}ce \`{a} nos r\'{e}gulateurs
intelligents, les correcteurs PI et PID classiques, \'{e}chantillonn\'{e}s,
donnent des r\'{e}sultats satisfaisants pour des syst\`{e}mes complexes si
leur r\'{e}glage est appropri\'{e}. Le b\'{e}n\'{e}fice, consid\'{e}rable, de nos iPI
est, comme d\'{e}j\`{a} dit, un r\'{e}glage simple.
\begin{remarque}
Voir \cite{pid} pour plus de d\'{e}tails. Rappelons que cette
compr\'{e}hension, nouvelle semble-t-il, des PID a \'{e}t\'{e} d\'{e}gag\'{e}e lors de
l'implantation du sans-mod\`{e}le \`{a} un v\'{e}hicule \'{e}lectrique \cite{cifa}.
\end{remarque}

\section{Mod\`{e}le restreint}\label{nlspring}
Soit le syst\`eme masse-ressort:
\begin{equation}
m \ddot{y} = -k_1 y - k_3y^3+ {\mathcal{F}}(\dot{y}) - d\dot{y} +u\label{mass}
\end{equation}
o\`u
\begin{itemize}
\item seule la masse $m = 0.5$ est connue;
\item $k_1=3,$ $k_3=10$, $d=5$ sont suppos\'{e}s mal connus, et on utilise
les estim\'{e}es: $\hat{k}_1=2$, $\hat{k}_3=7$, $\hat{d}=2.5$;
\item un frottement discontinue
$$ {\mathcal{F}}(\dot{y}) = \left\{
\begin{array}{c}
-0.3-0.4(\dot{y} + 0.25)^2  - d \dot{y} ~~~ \text{si}~~~ \dot{y}>0 \\
0.3+0.4(\dot{y} + 0.25)^2  - d \dot{y} ~~~~~\text{si} ~~~\dot{y}<0
\end{array}
\right. $$ Son expression math\'{e}matique, indispensable pour la
simulation, est inconnue du correcteur.
\end{itemize}
On utilise la platitude (\cite{flmr,mfrm,levine,hsr}) du mod\`{e}le
lin\'{e}aire restreint
\begin{equation}
m \ddot{y} = -\hat{k}_1 y - \hat{k}_3y^3 - \hat{d}\dot{y} +u
\label{mass11}
\end{equation}
pour d\'{e}terminer une commande nominale en boucle ouverte:
\begin{equation*}
u^{\star} =m \ddot{y}^{\star} +\hat{k}_1 y^{\star}+\hat{k}_3
(y^{\star})^3 + \hat{d}\dot{y}^{\star}  \label{mass12}
\end{equation*}
\begin{figure}[htp]
\centering
\includegraphics[width=9.05cm]{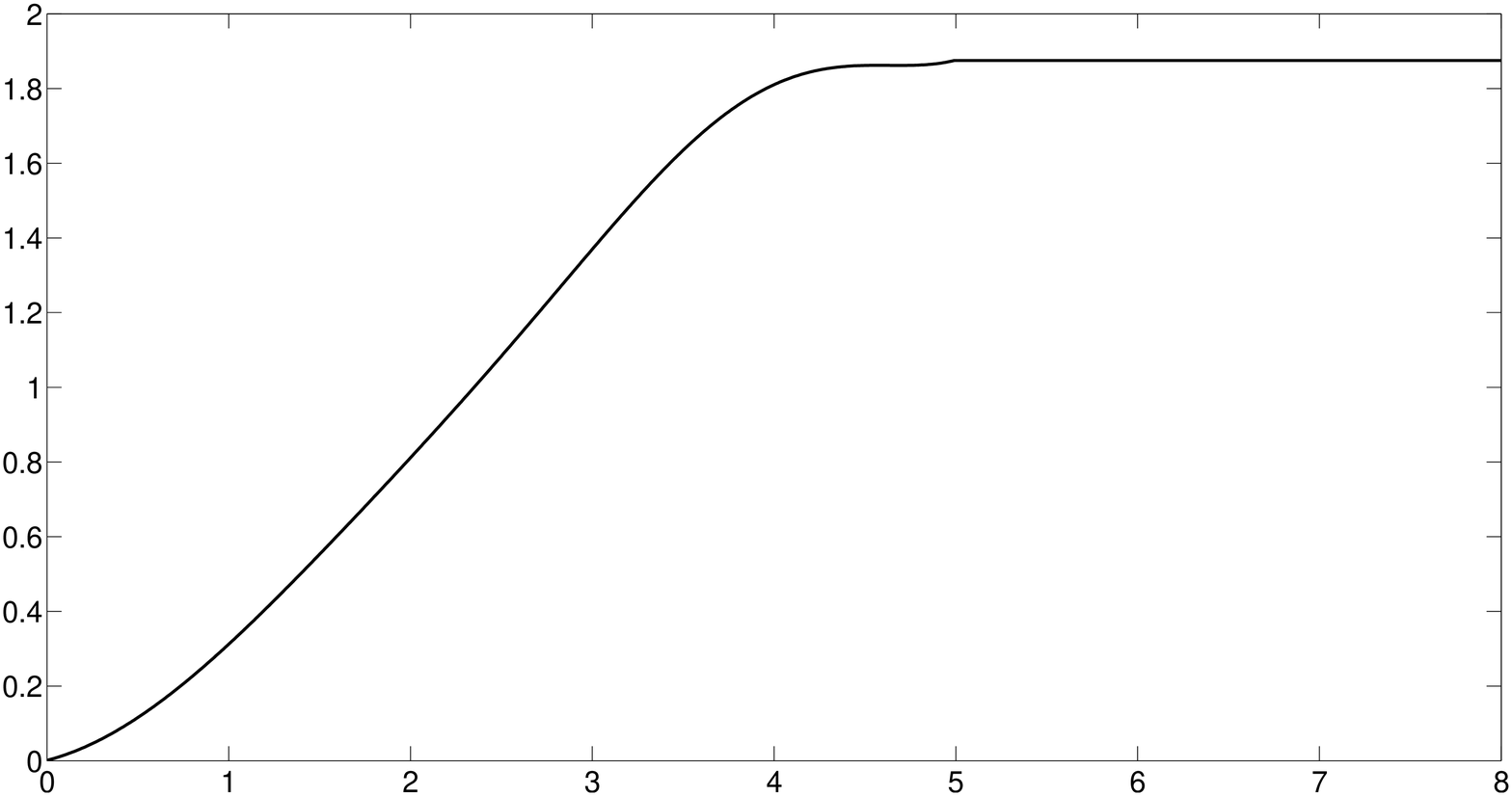}
\caption{Commande nominale du ressort.}\label{massa}
\end{figure}
\begin{figure}[htp]
\centering
\includegraphics[width=9.05cm]{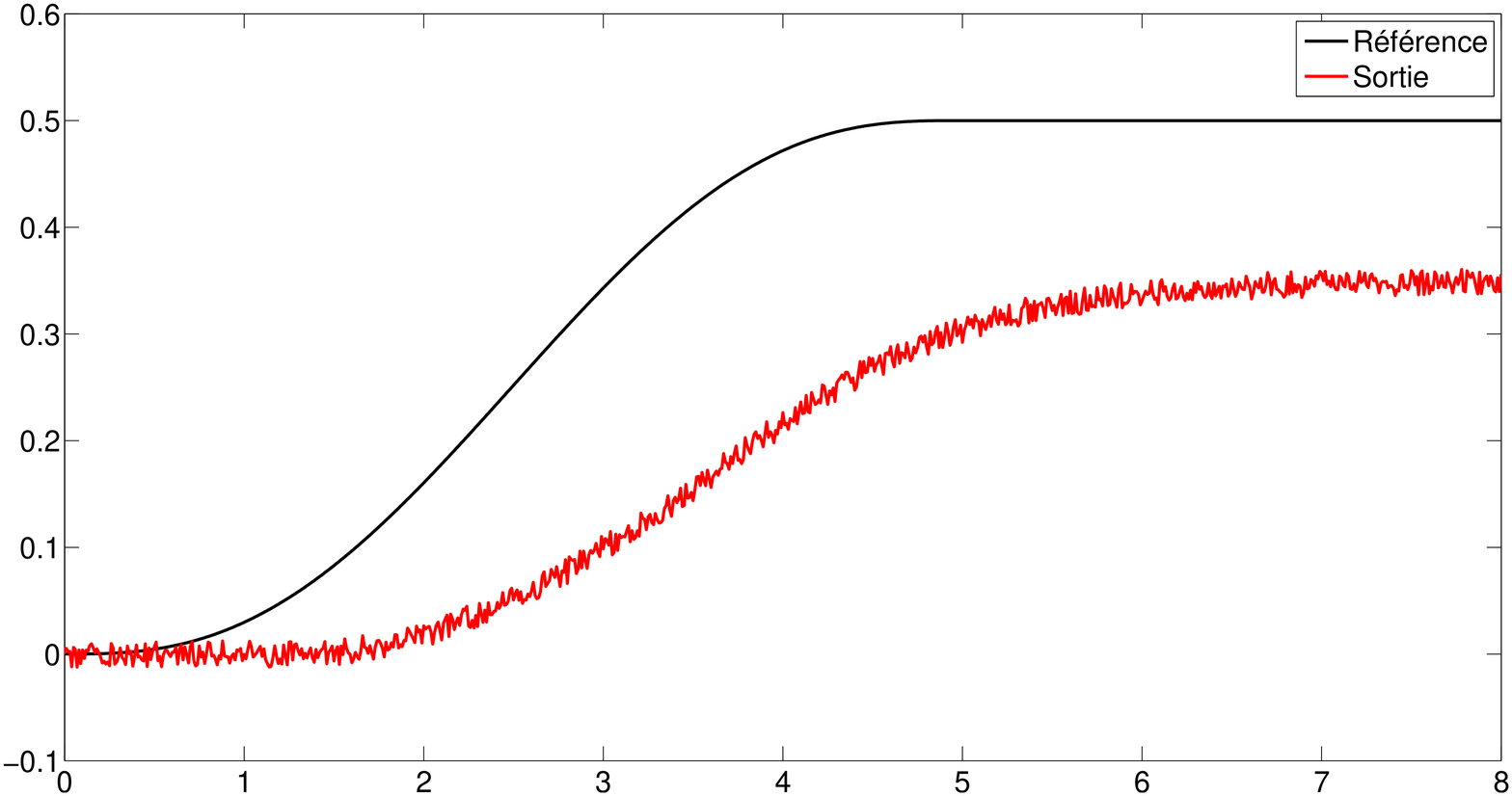}
\caption{Sortie avec commande nominale.}\label{massb}
\end{figure}
Les simulations des figures \ref{massa} et \ref{massb} correspondent
\`a $u^\star$. Posons $u = u^\star + \Delta u$ et introduisons un
iPI par rapport \`{a} $\Delta u$ pour stabiliser $e = y - y^\star$
autour de $0$. On aboutit aux excellentes simulations des figures
\ref{massc} et \ref{massd}, obtenues avec du bruit.

\begin{figure}[htp]
\centering
\includegraphics[width=9.05cm]{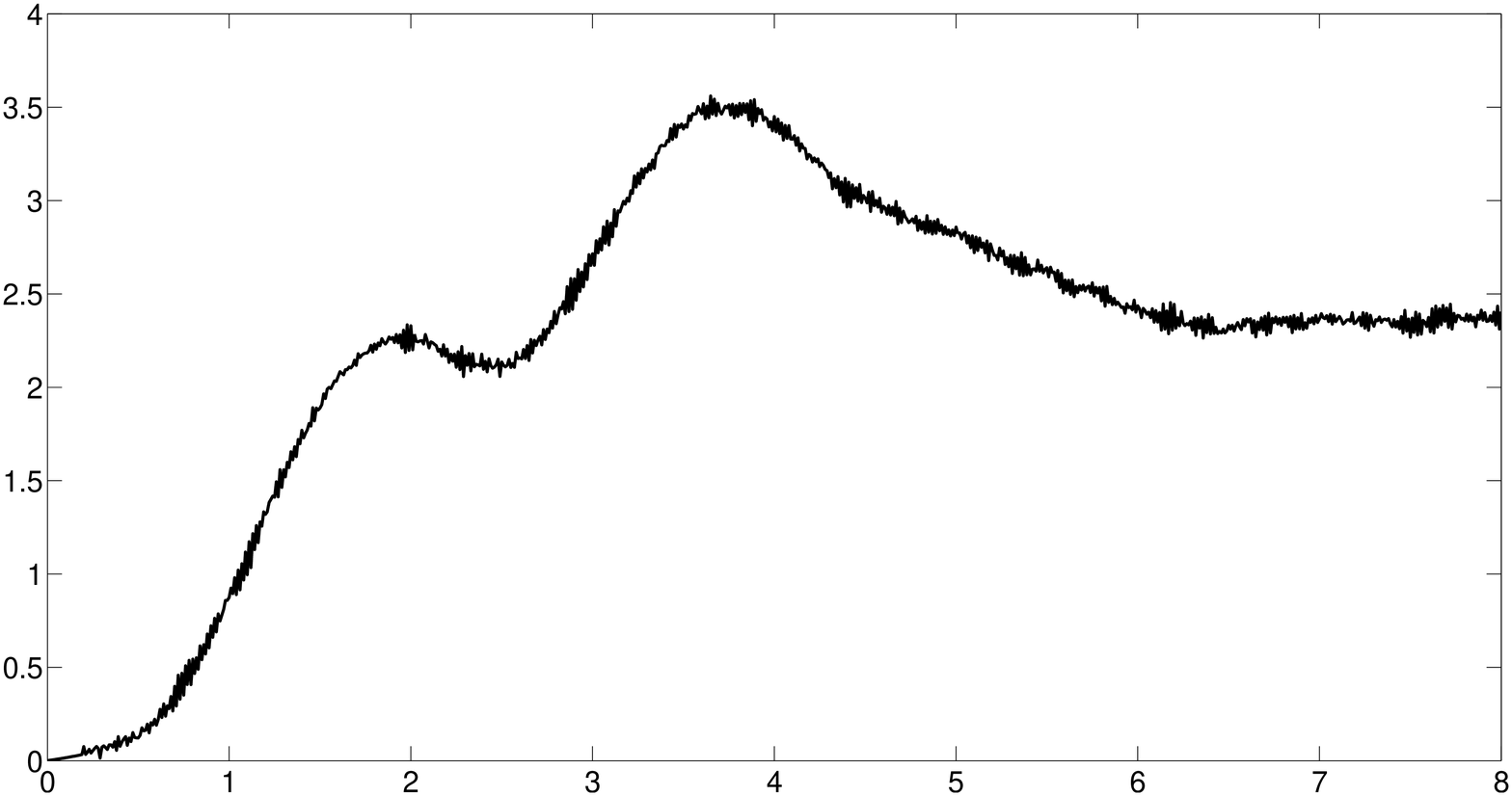}
\caption{Commande iPI du ressort dans le cas bruit\'e.}\label{massc}
\end{figure}

\begin{figure}[htp]
\centering
\includegraphics[width=9.05cm]{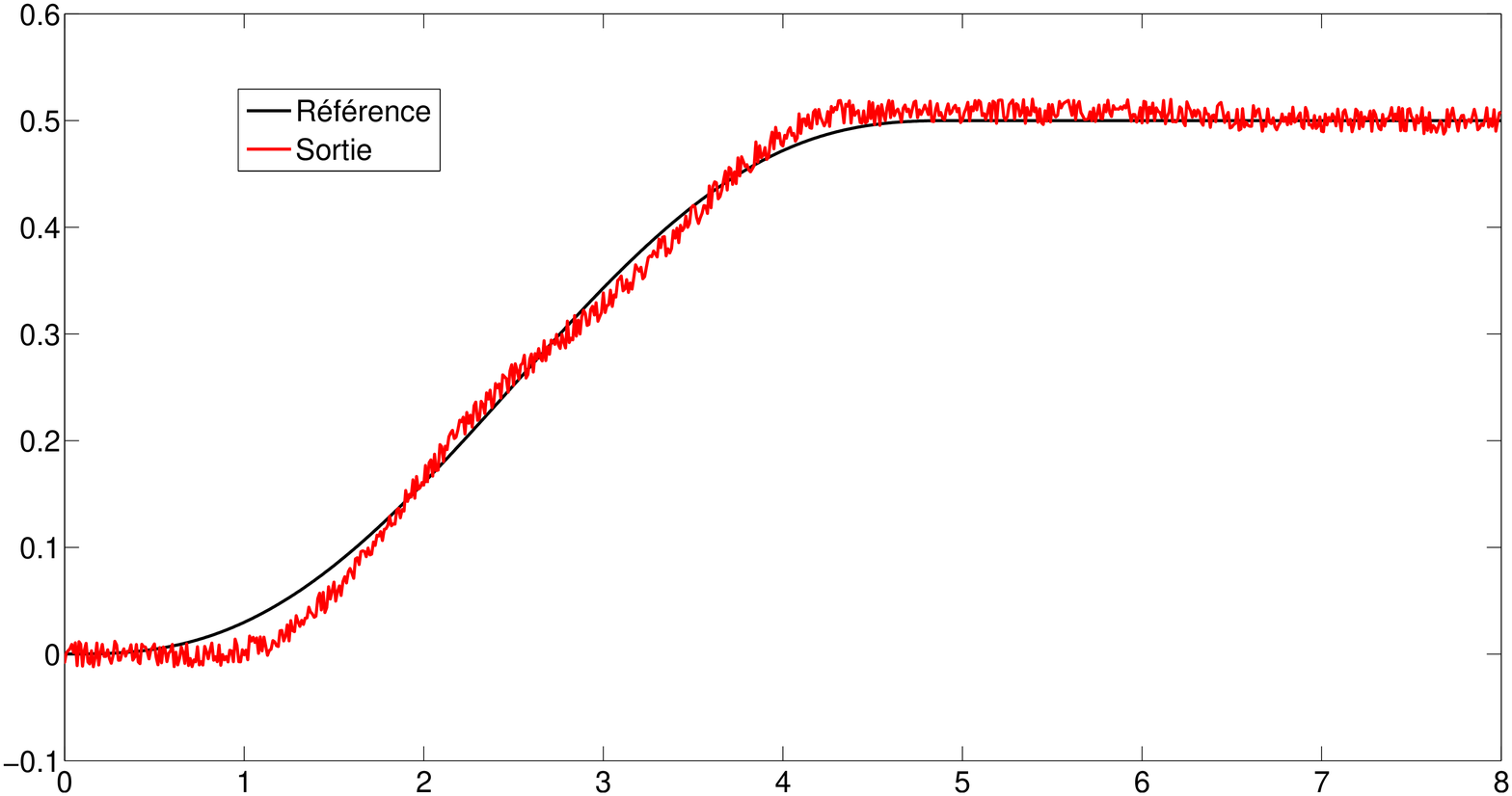}
\caption{La sortie du ressort dans le cas bruit\'e et commande
iPI.}\label{massd}
\end{figure}

\begin{remarque}
Le mod\`{e}le restreint \eqref{mass11} sert uniquement, comme on vient
de le voir, \`{a} calculer ais\'{e}ment une trajectoire de r\'{e}f\'{e}rence en
boucle ouverte, t\^{a}che ais\'{e} s'il est plat.
\end{remarque}

\section{Autres illustrations num\'{e}riques}\label{illu}

\subsection{R\'{e}gimes glissants}\label{autonome}
\begin{figure}[htp]
\centering
\includegraphics[width=9.05cm]{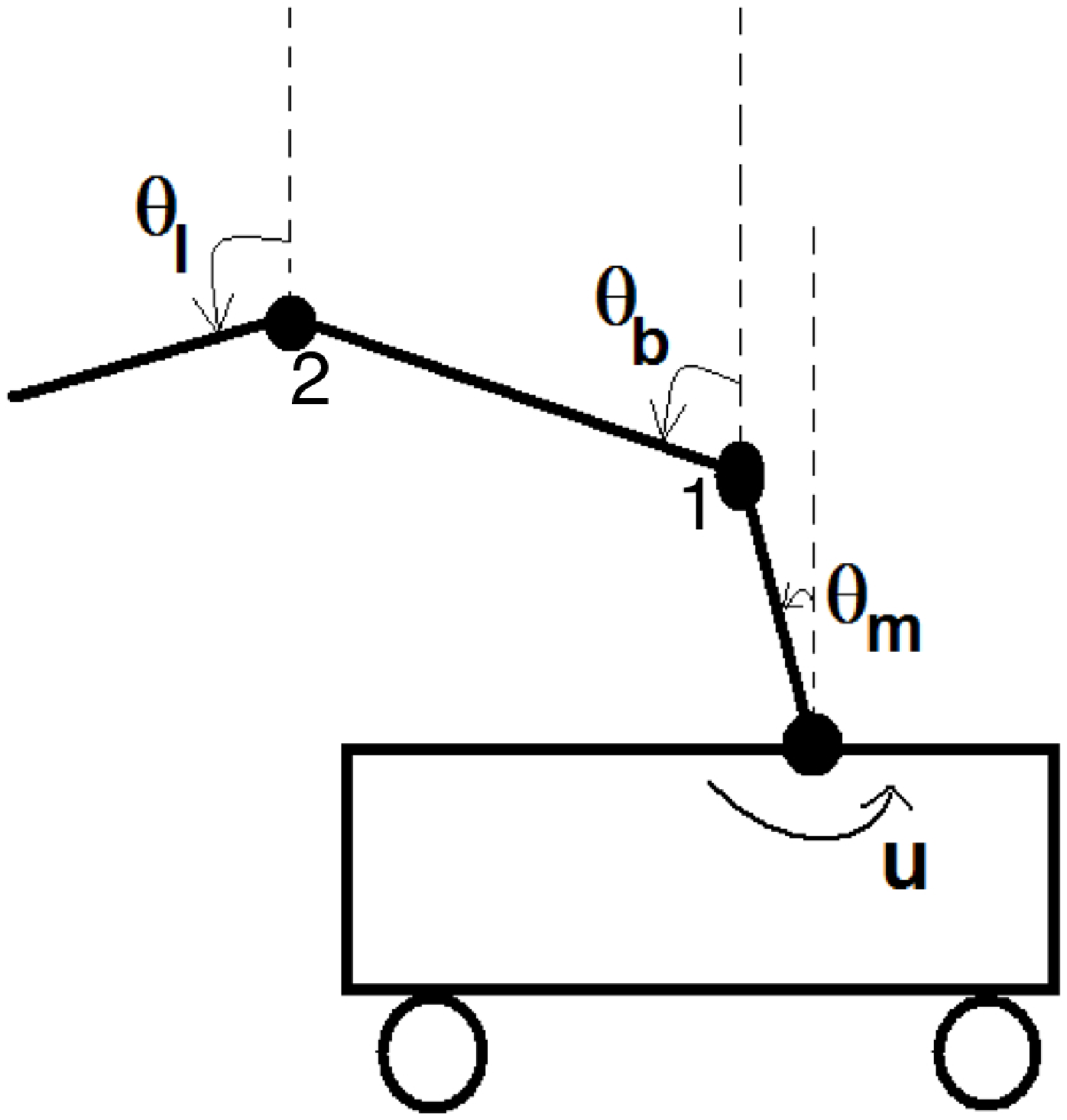}
\caption{V\'ehicule militaire autonome.}\label{dessin}
\end{figure}
La figure \ref{dessin} repr\'{e}sente un v\'{e}hicule militaire autonome,
emprunt\'{e} \`{a} \cite{spurgeon}. Son int\'{e}r\^{e}t est d\^{u} au non respect de la
\og \textit{matching condition} \fg, bien connue des tenants de la
commande par r\'{e}gimes glissants (voir, par exemple,
\cite{draz,leonid2,leonid,spurgeon,lille}), ce qui entra\^{\i}ne une mise
en place p\'{e}nible de cette commande en \cite{spurgeon}, alors que le
sans-mod\`{e}le s'applique sans difficult\'{e} aucune. Avec une dynamique
lin\'{e}aire du type
\begin{equation}
\dot{x} = A x + Bu + B_{1}\omega \label{Spurgeon2}
\end{equation}
o\`{u} $\omega$ est une perturbation pouvant d\'{e}pendre du temps,
rappelons que cette condition correspond \`a $B_1=\beta B$, o\`u
$\beta$ est une constante.

\begin{figure}[htp]
\centering
\includegraphics[width=9.05cm]{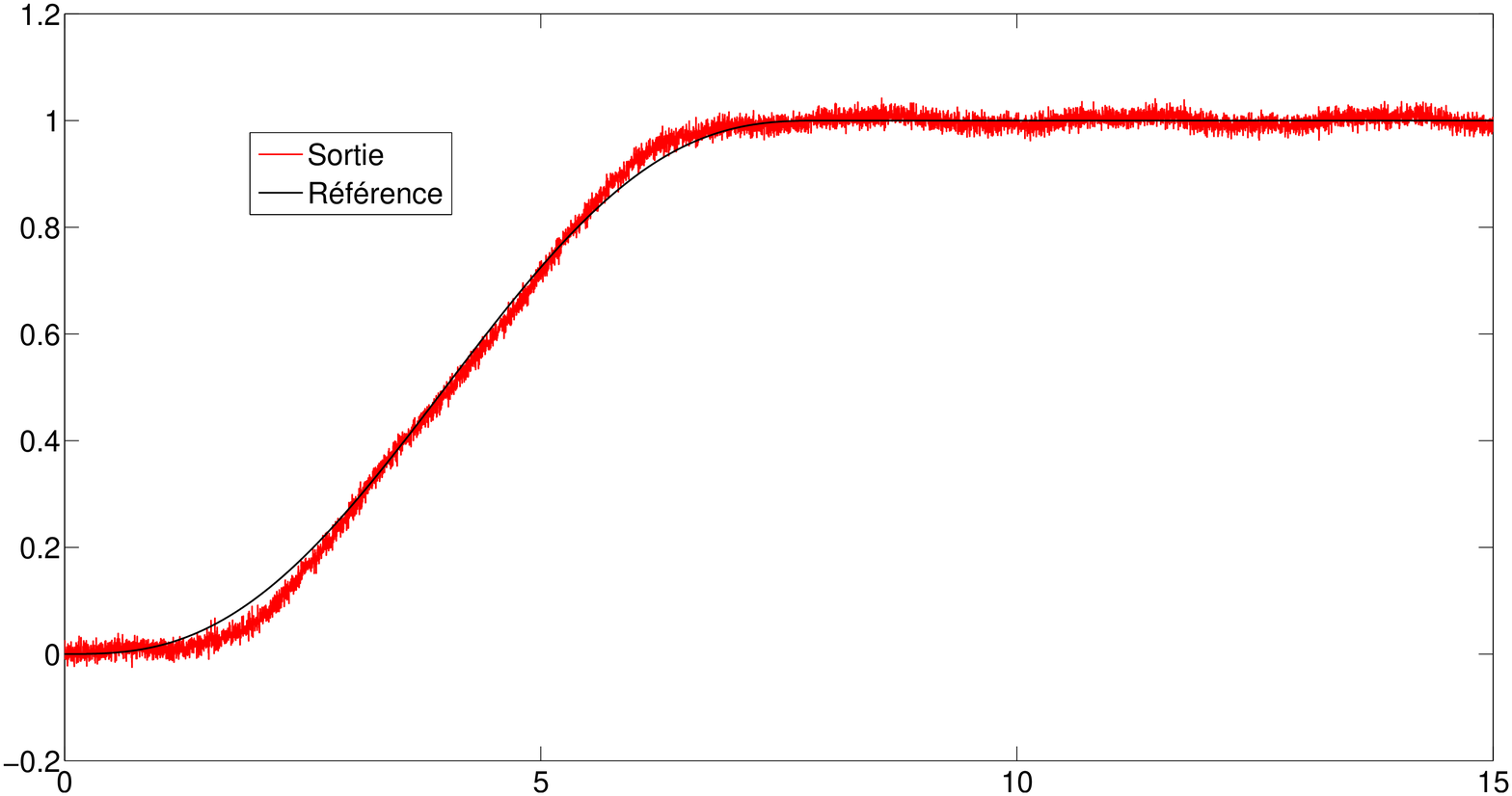}
\caption{Angle $\theta_b$.}\label{fig9}
\end{figure}

\begin{figure}[htp]
\centering
\includegraphics[width=9.05cm]{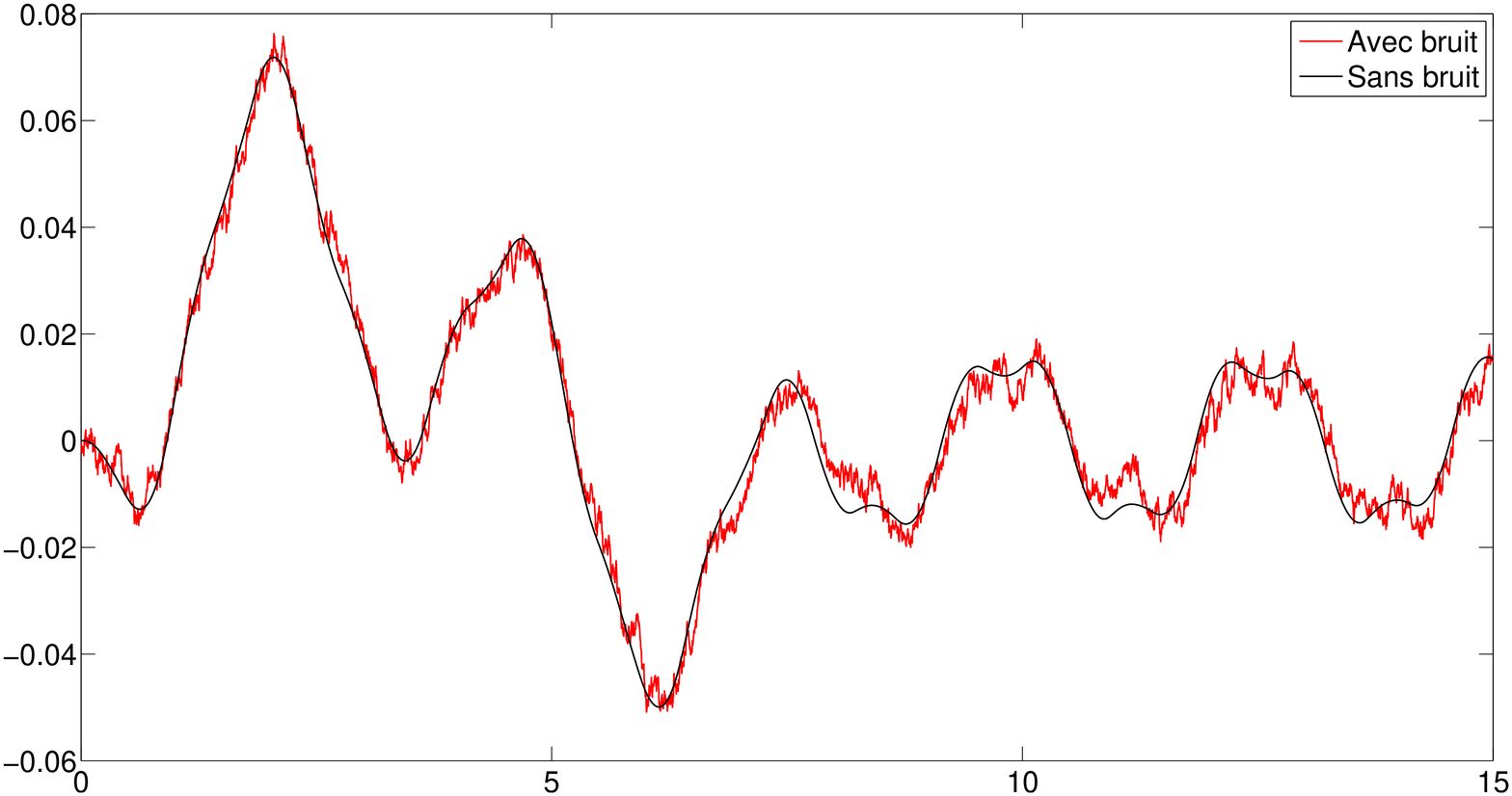}
\caption{Commande.}\label{fig10}
\end{figure}

L'angle $\theta_b$ doit prendre une valeur donn\'{e}e avec une commande
$u$ agissant sur $\theta_m$ et une perturbation issue $\theta_l$.
Les raideurs et les coefficients d'amortissement des
ressorts-amortisseurs en les articulations 1 et 2, donn\'{e}s par la
matrice $A$, sont assez \'elev\'ees. Les \'equations dynamiques sont
de la forme (\ref{Spurgeon2}):
$x=[\theta_b,\dot{\theta}_m,\theta_{mb},\dot{\theta}_b,\theta_{bl},\dot{\theta}_l]^T$,
$\omega(t) = [0,\dot{\theta}_p,\omega_{1m},\omega_{1l},F_{mb}\times
sign(\dot{\theta}_m-\dot{\theta}_p),F_{mb}\times
sign(\dot{\theta}_b-\dot{\theta}_p)]^T$, o\`{u}
\begin{itemize}
\item $\dot{\theta}_m$:  vitesse angulaire du moteur,
\item $\dot{\theta}_b$: vitesse angulaire du bras interm\'ediaire,
\item $\theta_p$: angle de tangage du v\'ehicule pris comme
perturbation,
\item
$\theta_{mb}=\frac{1}{N}(\theta_m-\theta_p)+\theta_p-\theta_b$,
\item $N$: rapport de r\'eduction du moto-r\'educteur,
\item $\dot{\theta}_l$: vitesse angulaire du troisi\`eme bras
$\theta_l$,
\item $\omega_{1m}$: frottement au niveau du moteur,
\item $\omega_{1l}$: frottement au niveau du troisi\`eme bras
$\theta_l$,
\item $\tau_{am}$: couple appliqu\'e au moteur,
\item $\tau_{al}$: couple appliqu\'e \`a la charge,
\item $f_d$: frottement,
\item $ \omega_{1m}(t) = f_d sign(\tau_{am}-J_m\ddot{\theta}_p)$,
\item $\omega_{1l}(t) = f_d sign(\tau_{al}-J_l\ddot{\theta}_p)$,
\item $J_m,$ $J_l$: inerties.
\end{itemize}
On garde les valeurs num\'eriques de \cite{spurgeon} pour $A,$ $B$,
$B_1$, $C$: {\small
\begin{equation} A = \left[ \begin{array}{cccccc}
0&0&0&1&0&0\\
0&-338.14&-2.55\times 10^7&50942&0&0\\
0&0.0066&0&-1&0&0\\
0&0.66&5\times 10^4&-110.1&-15\times 10^3&10\\
0&0&0&1&0&-1\\
0&0&0&7.69&11538&-7.69
\end{array}\right]\nonumber
\end{equation}}
\begin{equation}
B = \left[ \begin{array}{c}
0\\
4523.1\\
0\\
0\\
0\\
0
\end{array}\right]\nonumber
\end{equation}
\begin{equation}
B_1 = \left[ \begin{array}{cccccc}
0&0&0&0&0&0\\
0&-50604&-769&0&5.1&0\\
1&0&0&0&0&0\\
0&99&0&-0.01&0&0.01\\
0&0&0&0&0&0\\
0&0&0&0&0&0
\end{array}\right]\nonumber
\end{equation}
\begin{equation}
C = \left[ \begin{array}{cccccc}
1&0&0&0&0&0\end{array}\right]\nonumber
\end{equation}
Avec $\dot{\theta}_b = F + u$ en guise de \eqref{F}, un correcteur
iPI, o\`{u} $K_P=160$, $K_I=6400$, fournit de bons r\'esultats, d'apr\`{e}s
les figures \ref{fig9} et \ref{fig10}. Notre commande, fort simple,
est robuste par rapport aux bruits de mesure, non consid\'{e}r\'{e}s en
\cite{spurgeon}. Elle ne n\'{e}cessite pas, contrairement \`{a}
\cite{spurgeon},
\begin{itemize}
\item la r\'esolution d'in\'egalit\'es matricielles,
\item la mesure de $\dot{\theta}_m$, $\theta_{mb}$ et $\dot{\theta}_b$.
\end{itemize}

\subsection{Dimension infinie}
\subsubsection{Retards}\label{retards}
\begin{figure}[htp]
\centering
\includegraphics[width=9.05cm]{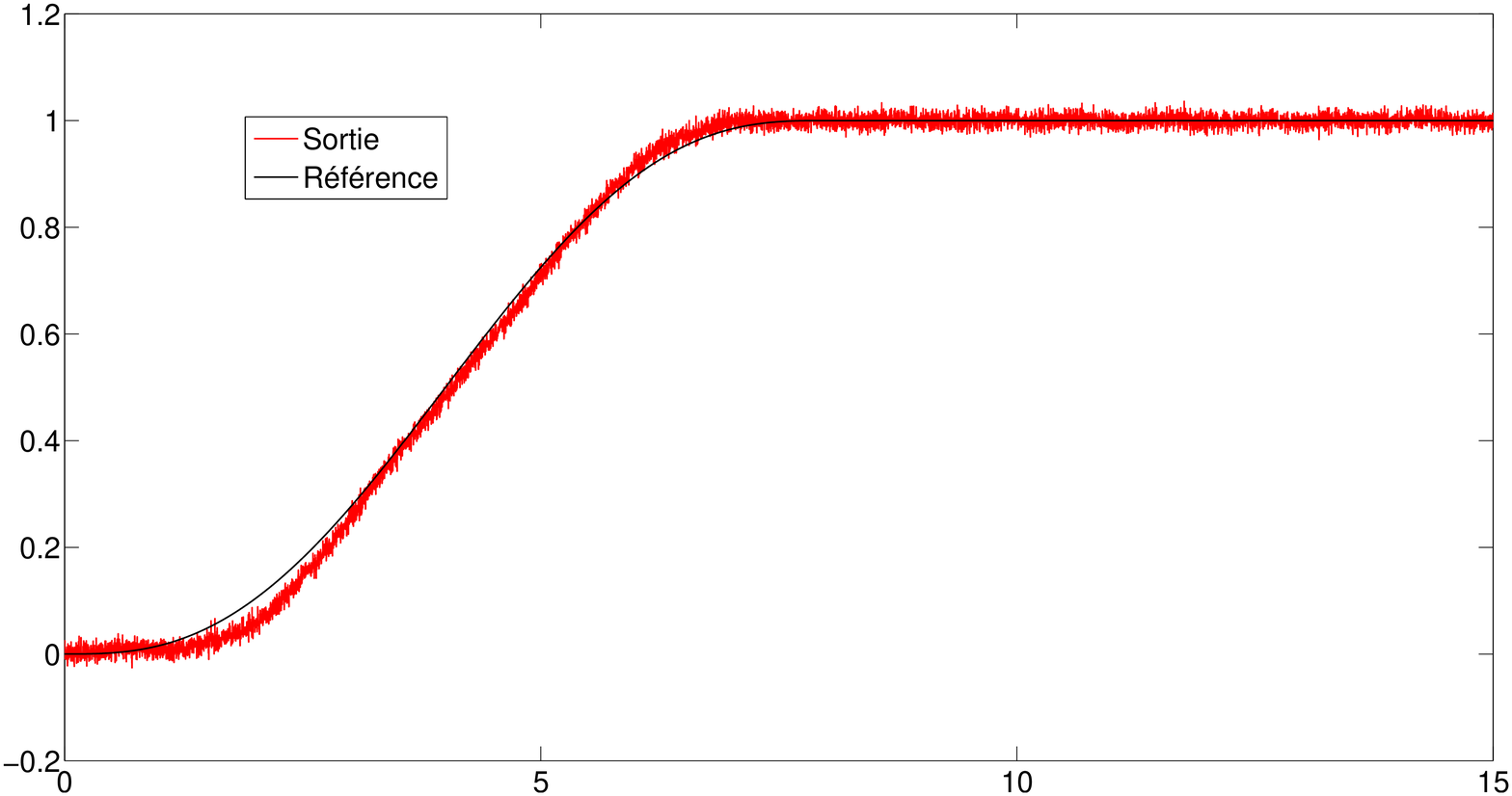}
\caption{Angle $\theta_b$.}\label{frid3}
\end{figure}

\begin{figure}[htp]
\centering
\includegraphics[width=9.05cm]{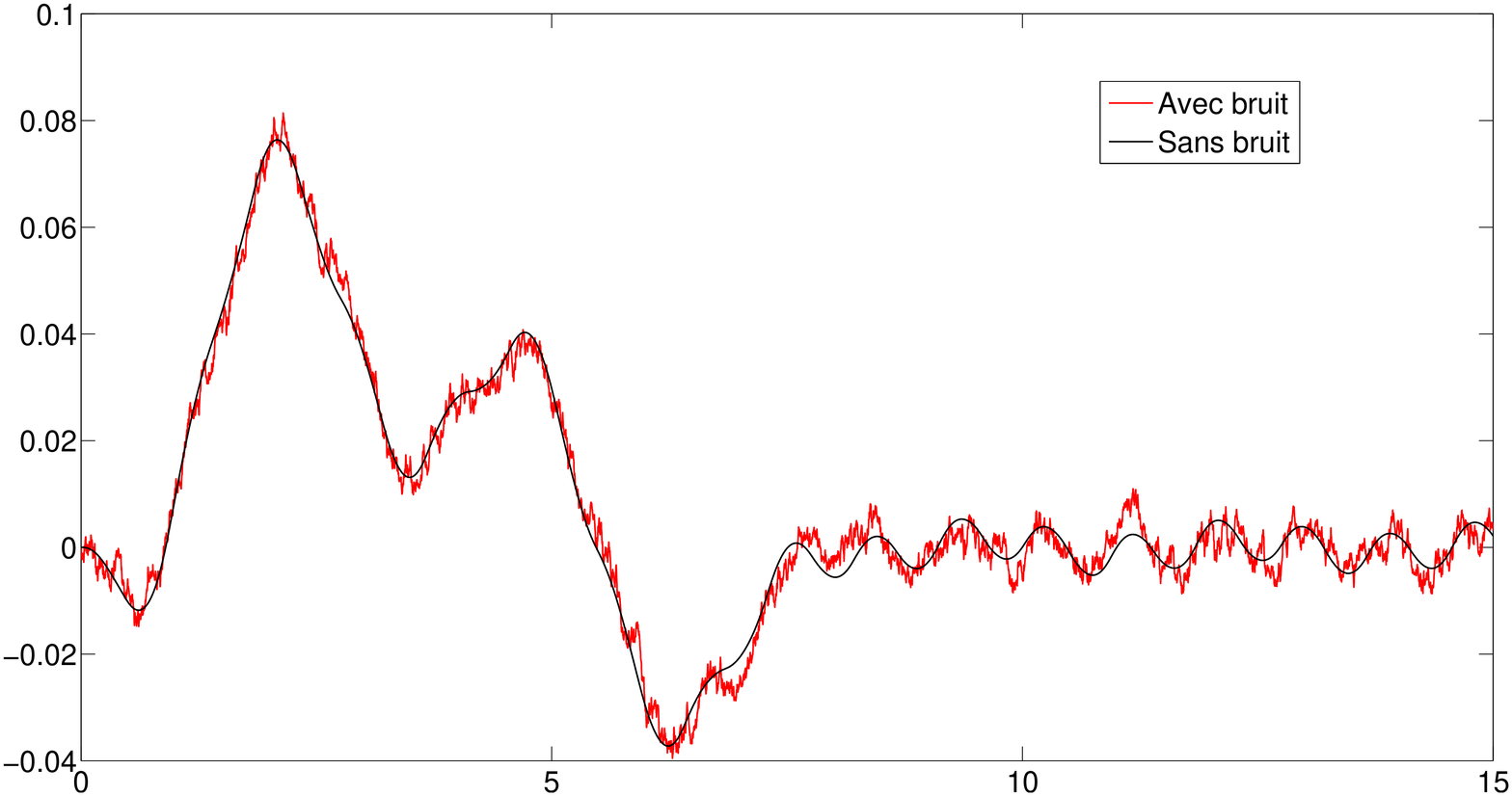}
\caption{Commande.}\label{frid4}
\end{figure}

Continuons l'exemple du {\S} \ref{autonome}. Ajoutons, comme
\cite{spurgeon}, un retard sur l'\'{e}tat, d\^{u} aux frottements. Alors,
\eqref{Spurgeon2} devient:
\begin{equation}\label{delay}
\dot{x}(t) = A x(t) + A_d x(t-\tau) + Bu(t) + B_{1}\omega(t)
\end{equation}
o\`u
\begin{itemize}
\item le retard $\tau$ varie entre $0$ et $9$ ms (l'amplitude est
moindre en \cite{spurgeon}, o\`{u} il varie entre $0$ et $3$ ms),
\item la matrice $A_d$, tir\'ee de \cite{spurgeon}, est
\begin{equation}
A_d = \left[ \begin{array}{cccccc}
0&0&0&0&0&0\\
0&-10&20&0&0&0\\
0&0.007&0&0.1&0&0\\
0&0&20&-2&1&0\\
0&0&0&0.1&0&0.1\\
0&0&0&2&1&0
\end{array}\right] \nonumber
\end{equation}
\end{itemize}
La m\^eme commande qu'au {\S} \ref{autonome}  fournit les r\'esultats
tout \`{a} fait satisfaisants des figures \ref{frid3} et \ref{frid4}.
Elle est, ici aussi, bien plus simple que celle propos\'{e}e en
\cite{spurgeon}.

\begin{remarque}
Les termes retard\'{e}s en \eqref{delay} se retrouvent en $F$ dans
\eqref{F}. Ils sont pris en compte sans autre forme de proc\`{e}s par
notre strat\'{e}gie de commande. Nul besoin donc d'identifier ces
retards, t\^{a}che d\'{e}licate, surtout s'ils sont variables.
\end{remarque}

\begin{remarque}
Rappelons (voir, par exemple, \cite{astrom,od}) que le r\'{e}glage des
PID traditionnels repose souvent sur une approximation du syst\`{e}me \`{a}
commander par un syst\`{e}me \`{a} retard. Ici aussi, ces retards sont
inutiles, comme il est d\'{e}montr\'{e} en \cite{esta,malo}, avec nos
correcteurs intelligents
\end{remarque}

\subsubsection{Syst\`emes \`a param\`etres distribu\'es}

\paragraph{\it Un \'{e}changeur de chaleur}

\begin{figure}[htp]
\centering
\includegraphics[width=9.05cm]{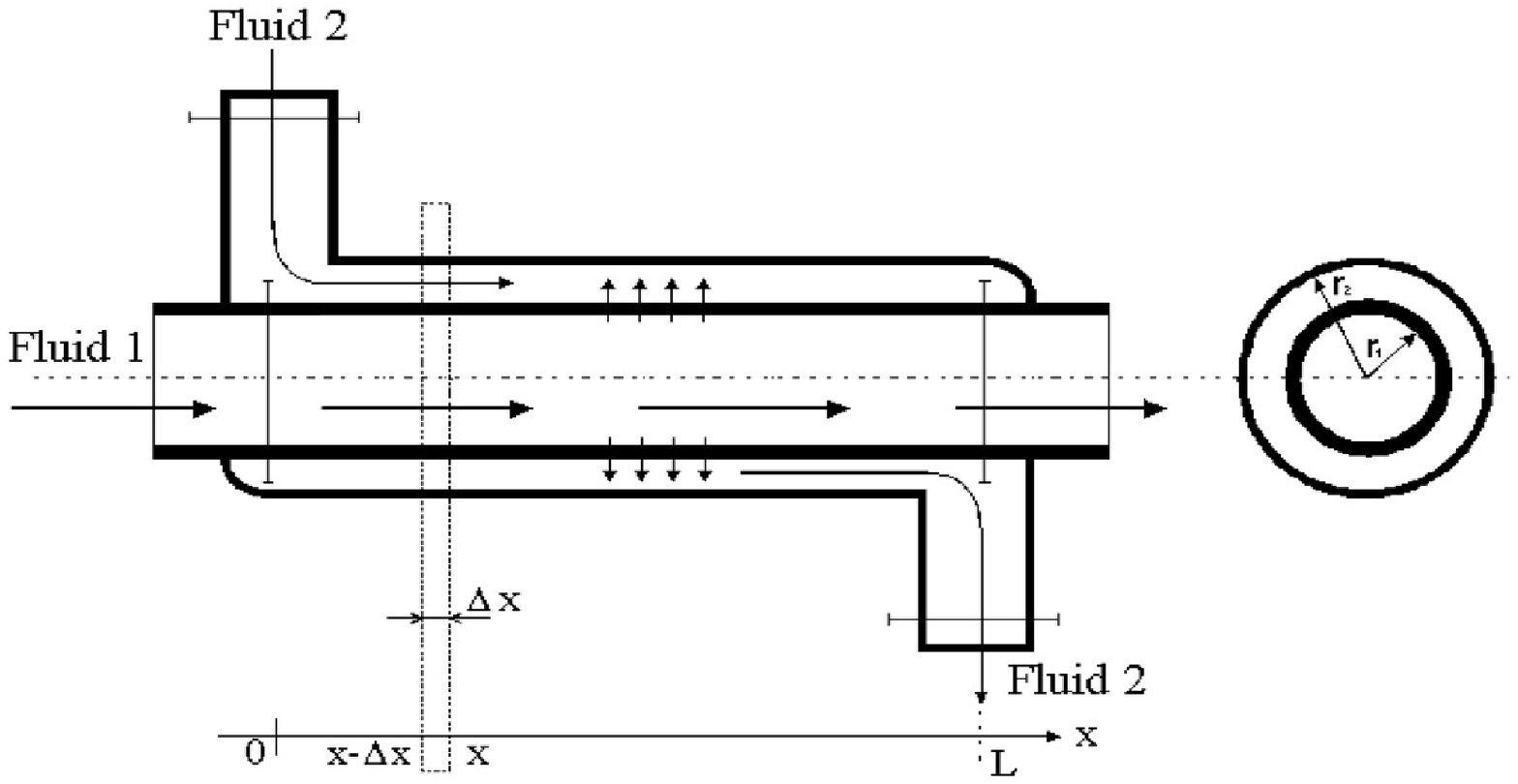}
\caption{\'Echangeur de chaleur.}\label{heat}
\end{figure}
Pour l'\'{e}changeur de chaleur de la figure \ref{heat}, nos simulations
num\'{e}riques utilisent une discr\'{e}tisation spatiale des \'{e}quations aux
d\'{e}riv\'{e}es partielles
\begin{eqnarray*}
M_1c_1\frac{\partial T}{\partial t} & = & m_1c_1\frac{\partial
T}{\partial
x}-AU(T-S)\nonumber \\
M_2c_2\frac{\partial S}{\partial t} & = & m_2c_2\frac{\partial
S}{\partial x} -AU (T-S),\label{equa1}
\end{eqnarray*}
o\`{u}
\begin{itemize}
\item $0 \leq x \leq L$;
\item $T(x,t)$ est la temp\'{e}rature du fluide $1$, et $T(0,t)$ la commande;
\item $S(x,t)$ est la temp\'erature du fluide $2$, et $S(L,t)$ la
sortie;
\item $U$ est la conduction.
\end{itemize}
Les valeurs num\'{e}riques des param\`etres sont emprunt\'ees \`a un
\'echangeur r\'eel \cite{echangeur5}.

\begin{figure}[htp]
\centering
\includegraphics[width=9.05cm]{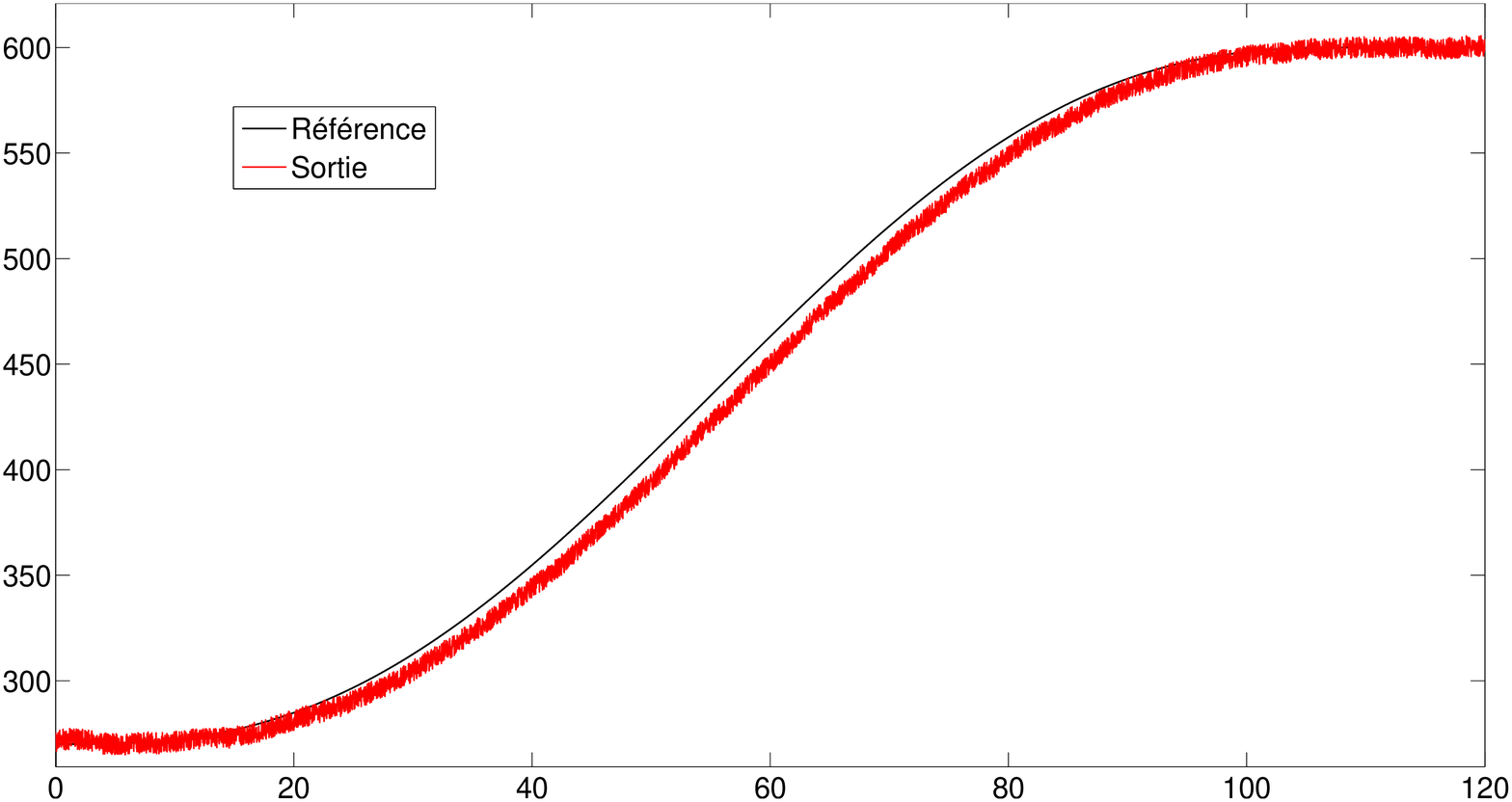}
\caption{$S(L,t)$ avec iPI.}\label{fig1}
\end{figure}

\begin{figure}[htp]
\centering
\includegraphics[width=9.05cm]{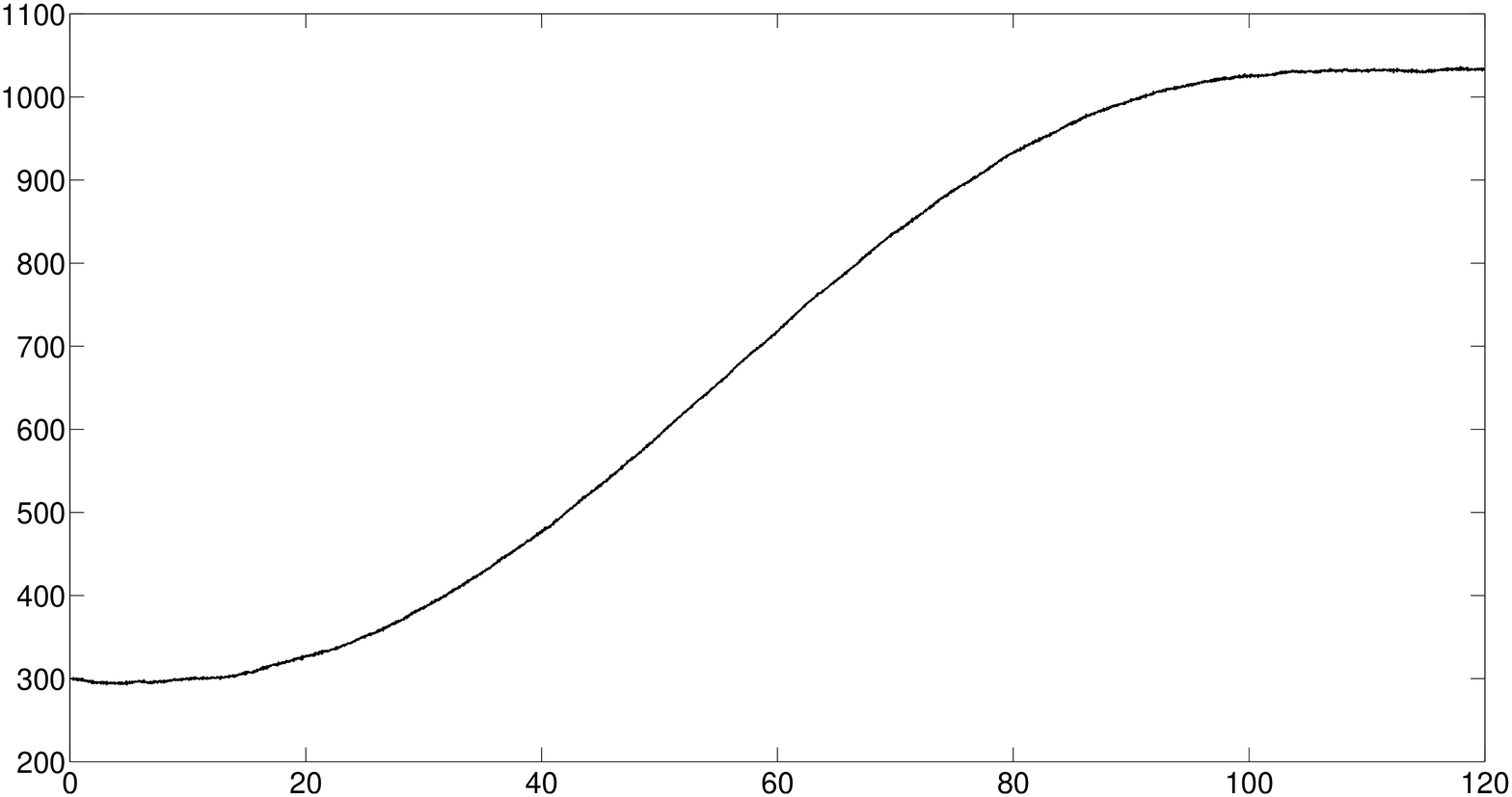}
\caption{Commande iPI.}\label{fig2}
\end{figure}

\begin{figure}[htp]
\centering
\includegraphics[width=9.05cm]{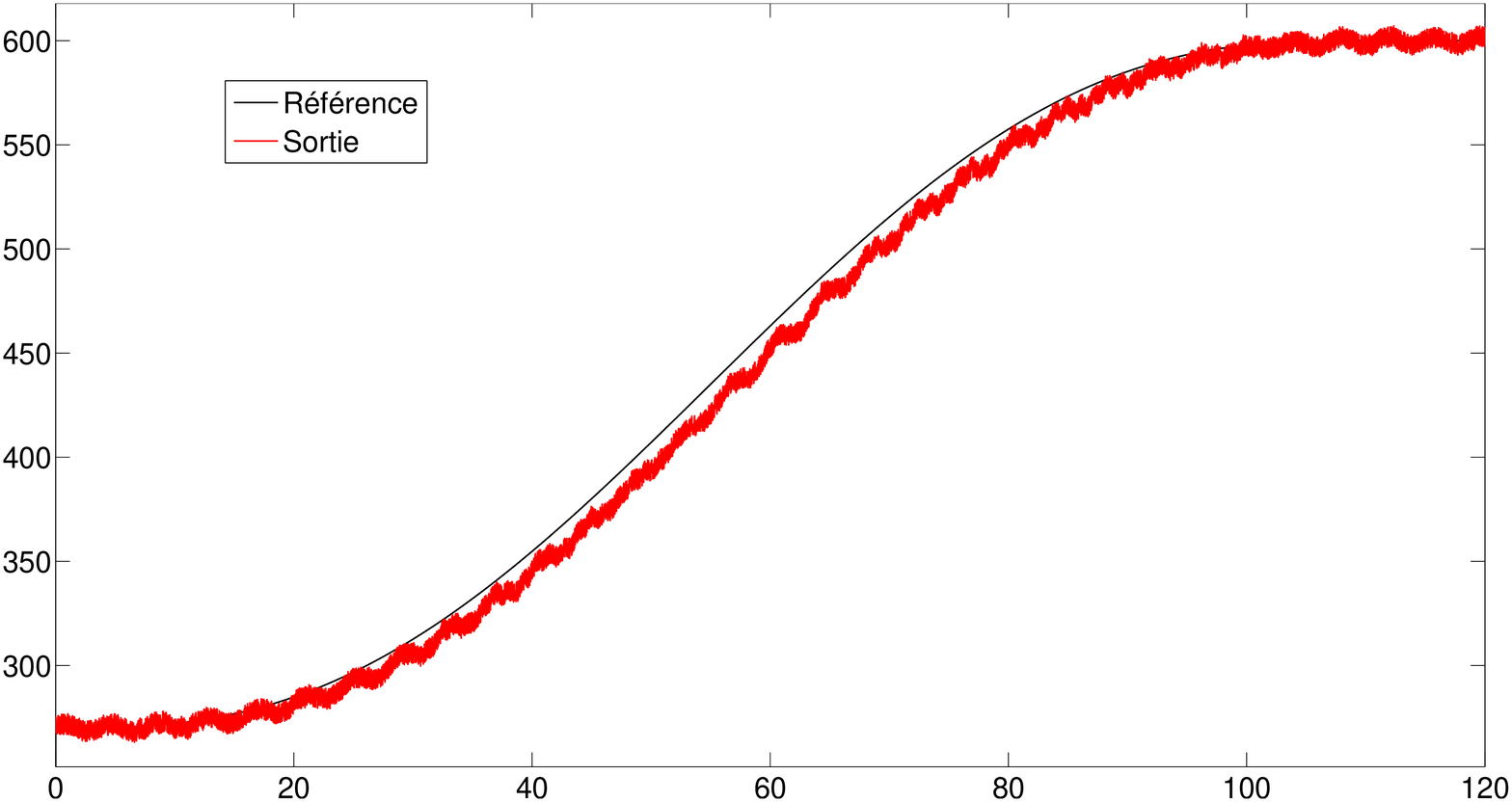}
\caption{$S(L,t)$ avec PID.}\label{fig3}
\end{figure}

\begin{figure}[htp]
\centering
\includegraphics[width=9.05cm]{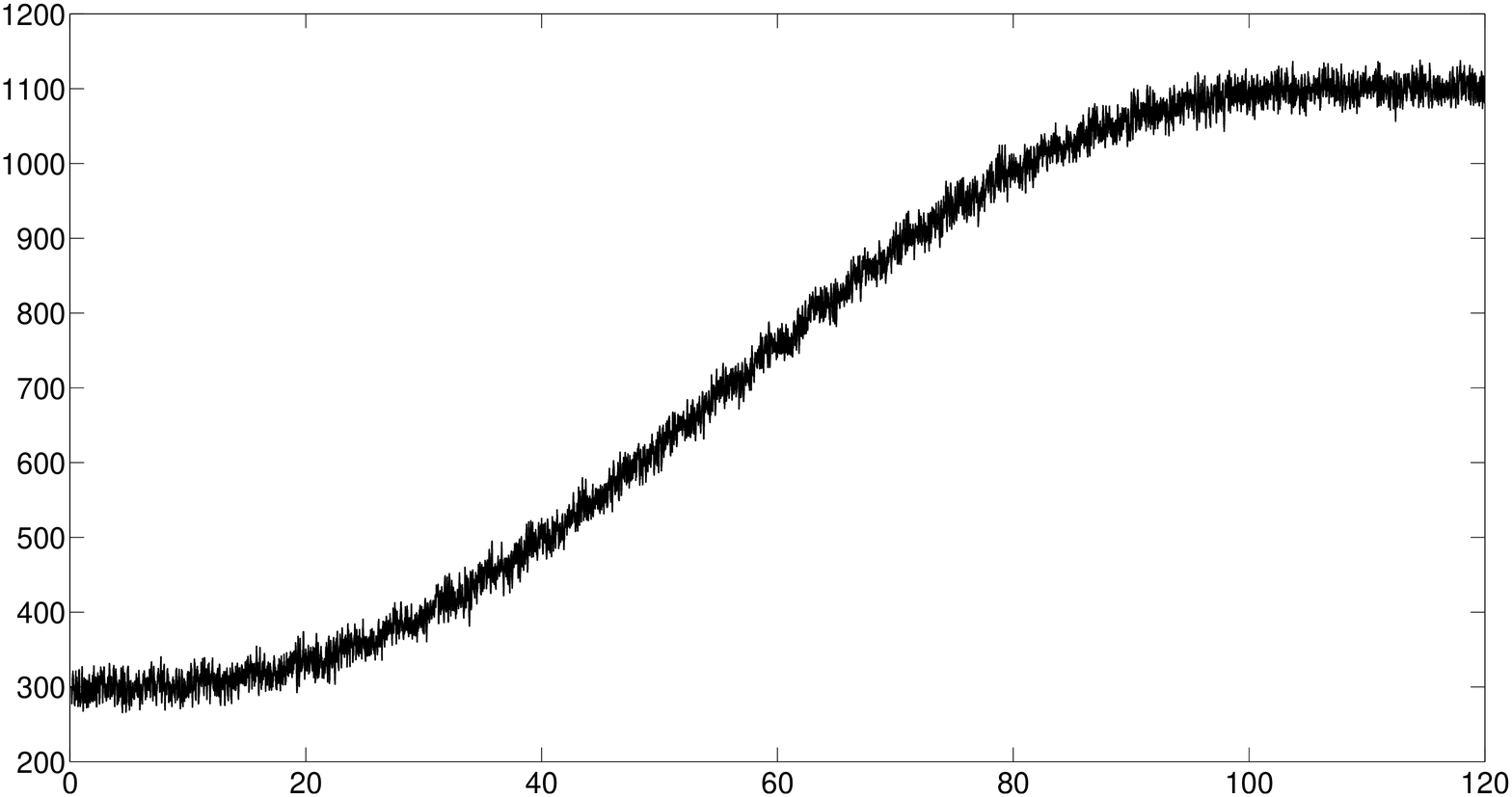}
\caption{Commande PID.}\label{fig4}
\end{figure}

\paragraph{\it iPI et PID} \eqref{F} devient:
\begin{equation}\label{18}
\dot{S}(L,t)=F+18T(0,t)
\end{equation}
On impose \`{a} $S(L,t)$ une trajectoire de r\'{e}f\'{e}rence douce allant de
$270^{\circ}$ \`a $600^{\circ}$ Kelvin, en utilisant une fonction polynomiale du
temps de degr\'{e} $6$. On \'evite ainsi l'excitation des modes rapides.
Quelques t\^{a}tonnements ont suffi pour arriver \`{a} la valeur $18$ en
\eqref{18}. Les gains de l'iPI, d\'{e}duit de \eqref{18}, sont $K_P=10$,
$K_I=17.36$. La m\'{e}thode de Ziegler-Nichols, bien connue et assez
lourde, fournit les trois gains du PID classique, construit pour les
besoins de la comparaison:
$$
T(0,t) = \dot{y}^{\star} + 1.8 \ e + \int e + 0.75 \ \dot{e}
$$
o\`{u} $y^{\star}$ est la trajectoire de r\'{e}f\'{e}rence et $e = y -
y^{\star}$. Les deux commandes fournissent de bons r\'{e}sultats comme
le montrent les figures \ref{fig1}, \ref{fig2}, \ref{fig3} et
\ref{fig4}. Notons, n\'{e}anmoins, la plus grande sensibilit\'{e} au bruit
du PID, due \`{a} la n\'{e}cessit\'{e} de calculer une d\'eriv\'ee.

\begin{figure}[htp]
\centering
\includegraphics[width=9.05cm]{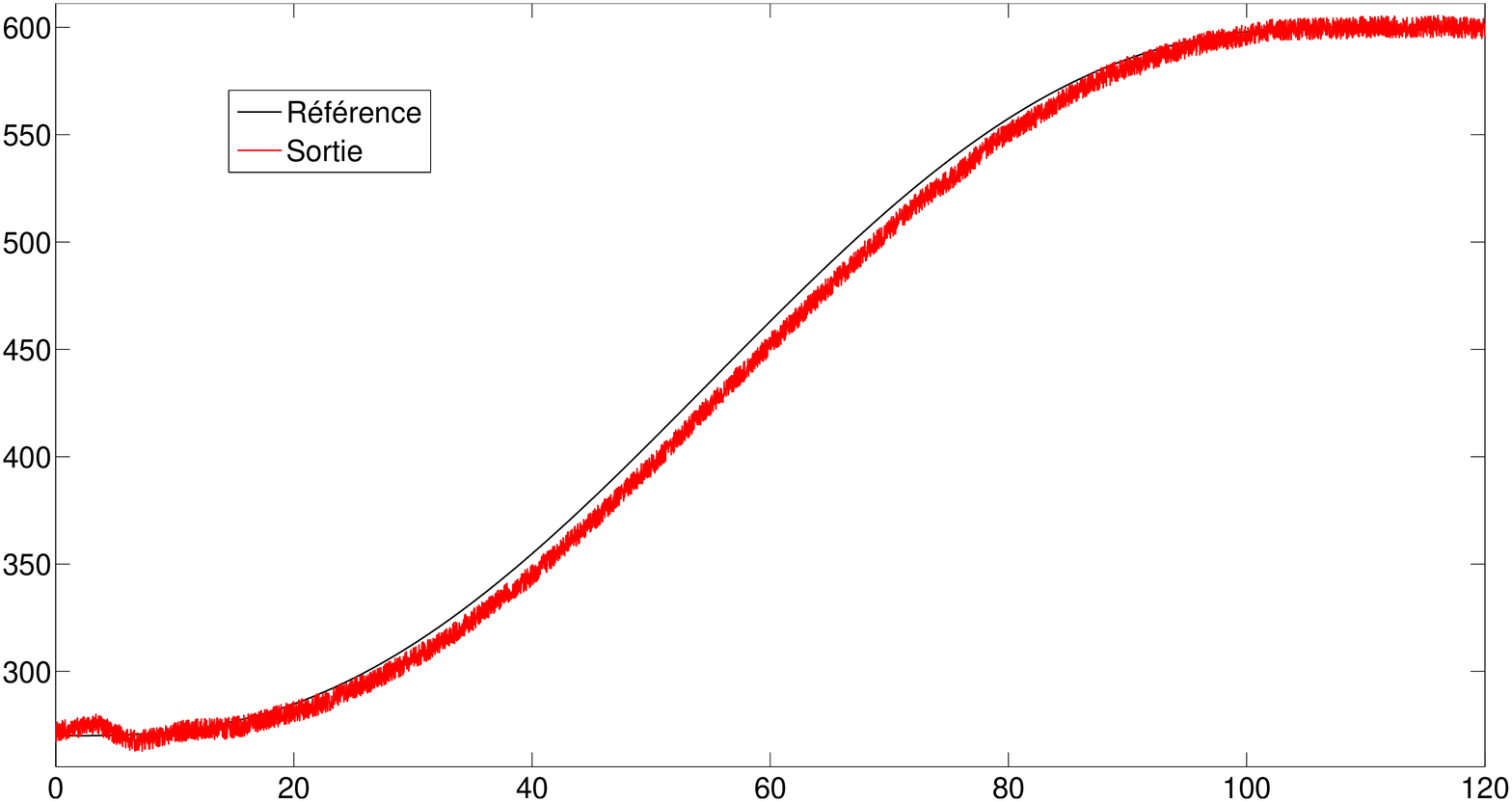}
\caption{$S(L,t)$ avec iPI.}\label{fig5}
\end{figure}

\begin{figure}[htp]
\centering
\includegraphics[width=9.05cm]{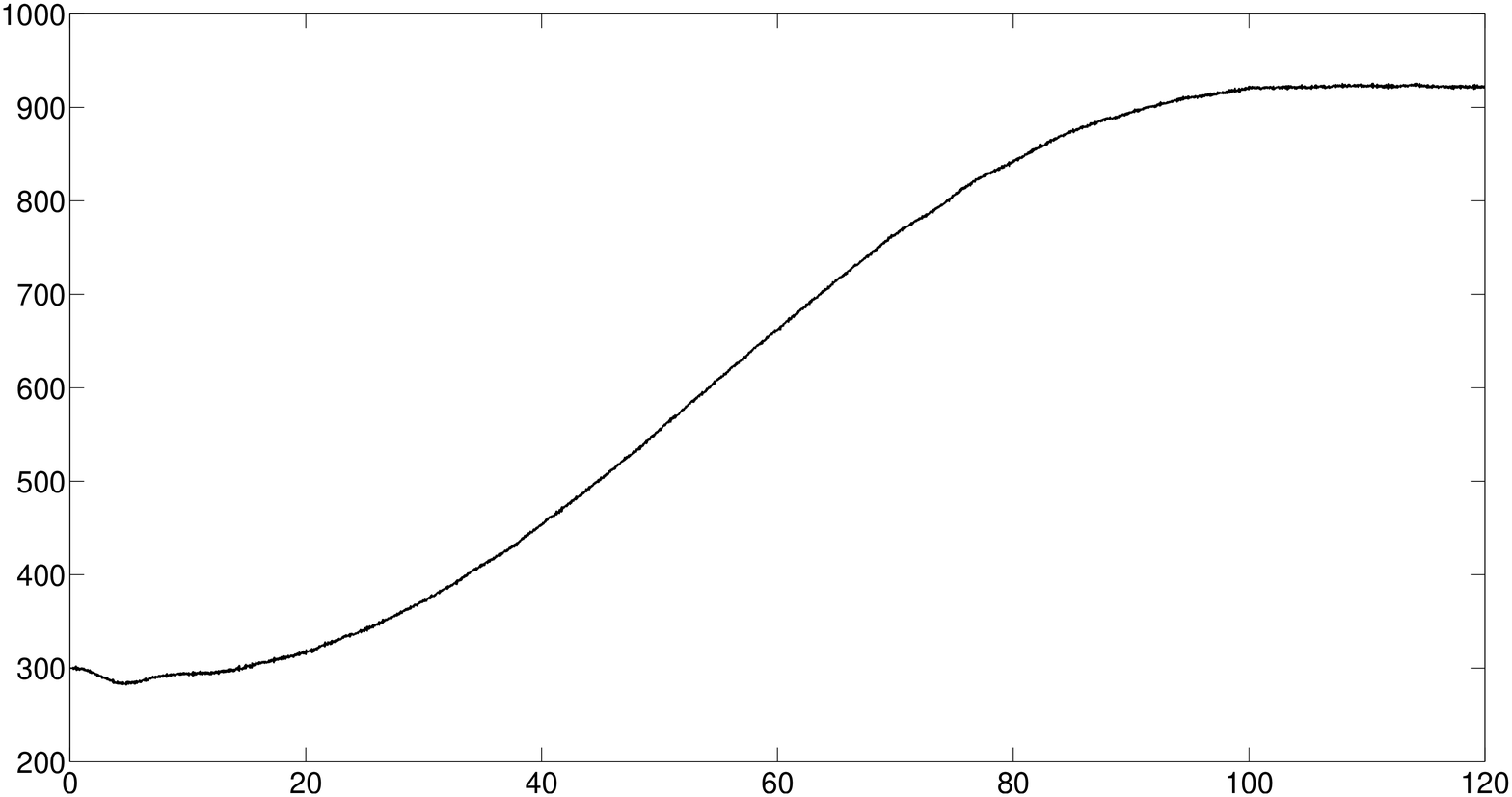}
\caption{Commande iPI.}\label{fig6}
\end{figure}

\begin{figure}[htp]
\centering
\includegraphics[width=9.05cm]{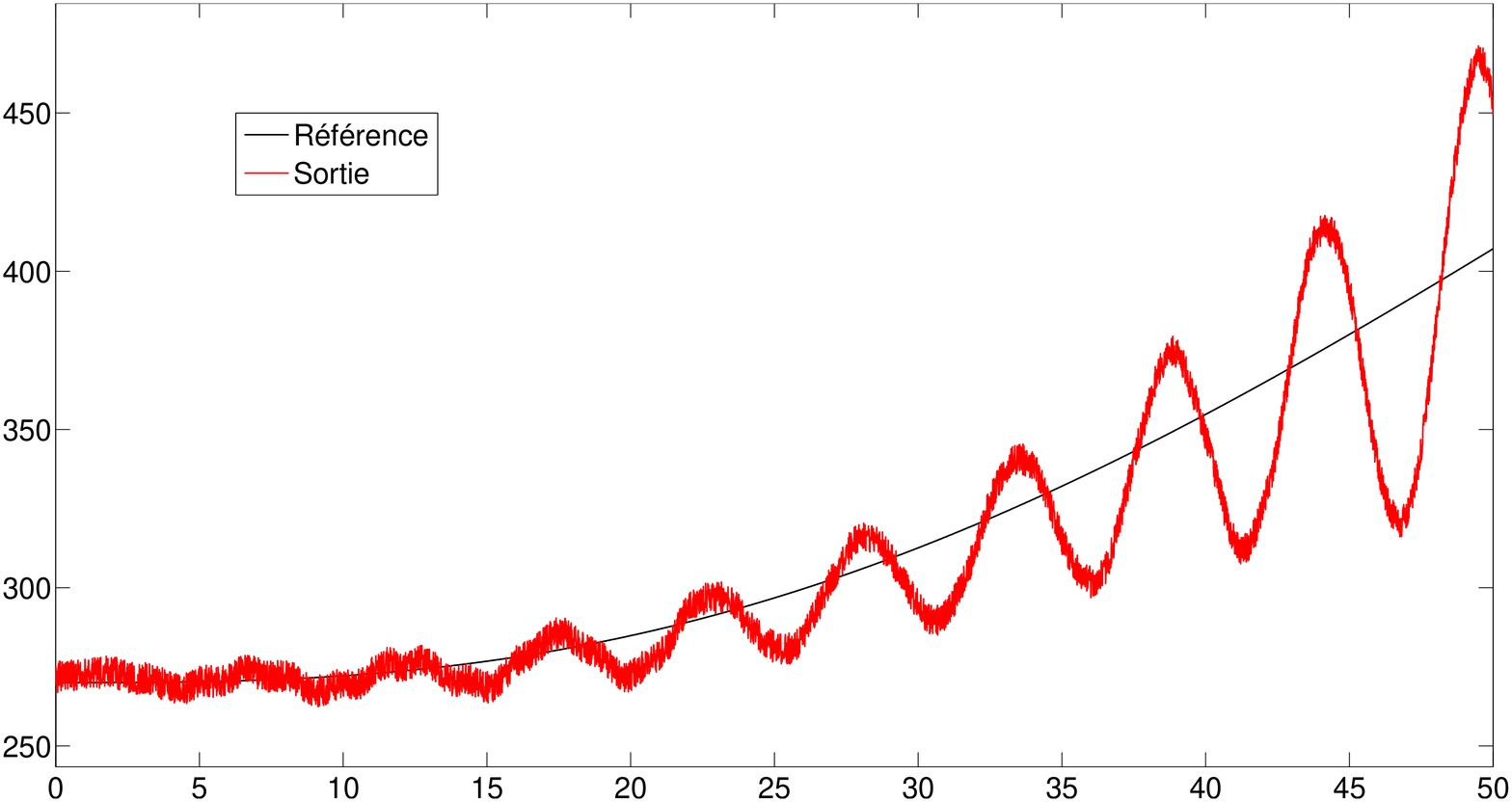}
\caption{$S(L,t)$ avec PID.}\label{fig7}
\end{figure}

\begin{figure}[htp]
\centering
\includegraphics[width=9.05cm]{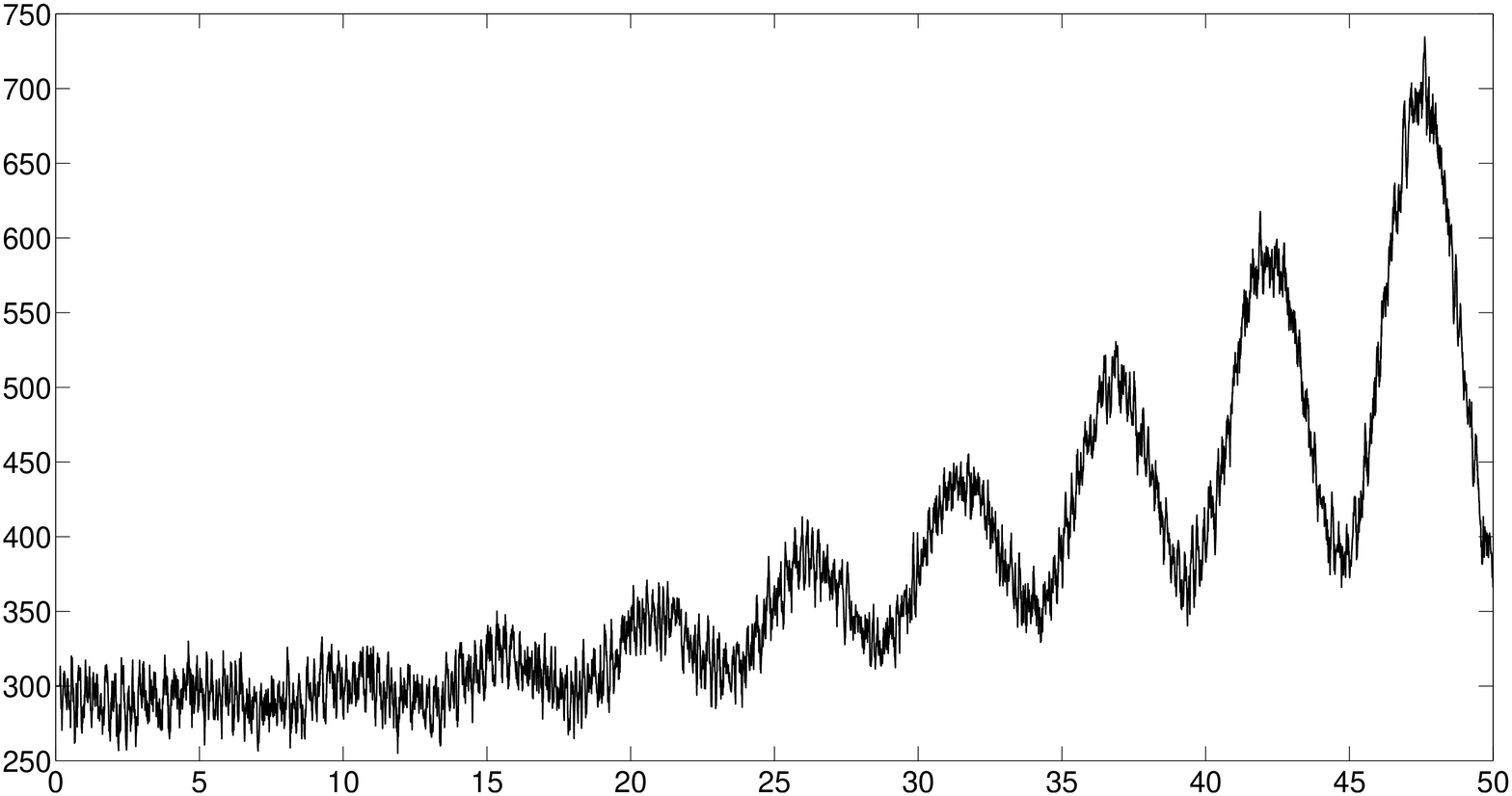}
\caption{Commande PID.}\label{fig8}
\end{figure}

Certains ph\'enom\`enes complexes sont in\'{e}vitables:
\begin{itemize}
\item les  vitesses et les propri\'et\'es physico-chimiques des fluides sont
variables;
\item le vieillissement s'accompagne de rouilles et de d\'{e}p\^{o}ts qui
peuvent modifier la conduction $U$ et r\'{e}tr\'{e}cir les tuyaux.
\end{itemize}
Toute commande de qualit\'{e} se doit de les ma\^{\i}triser. Supposons, donc,
les tuyaux partiellement bouch\'es, avec des baisses de $30\%$ du
d\'ebit et de $50\%$ de la conduction. Sans nouveau calibrage des
correcteurs pr\'{e}c\'{e}dents,
\begin{itemize}
\item les figures \ref{fig5} et \ref{fig6} montrent les performances
inchang\'ees de l'iPI;

\item celles du PID se d\'{e}gradent notablement d'apr\`{e}s les
figures \ref{fig7} et \ref{fig8}.
\end{itemize}

\paragraph{\it Autres commandes} Il y a une litt\'{e}rature importante sur la
commande des \'{e}changeurs de chaleur, qui ob\'{e}it \`{a} des points de vue
tr\`{e}s vari\'{e}s: voir\footnote{La liste ci-dessous n'a aucune pr\'{e}tention
d'exhaustivit\'{e}.}, par exemple,
\cite{echangeur5,echangeur1,heatexchbib,fisch,echangeur3,rudolph,echangeur2,edp2,echangeur4}.
Les \'{e}quations aux d\'{e}riv\'{e}es partielles, qui s'imposent dans toute
mod\'{e}lisation un tant soit peu fine, y sont \'{e}videmment pr\'{e}sentes,
mais ne semblent pas avoir produit une r\'{e}gulation exploitable.
Apparaissent aussi la commande pr\'{e}dictive, la commande robuste, la
lin\'{e}arisation par bouclage, la logique floue, \dots. Notre
correcteur iPI, simple \`{a} mettre en {\oe}uvre et aux bonnes performances,
semble plus prometteur.

\begin{remarque}
On trouve en \cite{horn} une int\'{e}ressante application du sans-mod\`{e}le
au thermique, qui souffre d'oscillations, dues, \`{a} notre avis, \`{a} une
d\'{e}rivation num\'{e}rique impropre. On y rem\'{e}die par les m\'{e}thodes du {\S}
\ref{oe}, par exemple.
\end{remarque}

\section{Conclusion}\label{conclusion}
Des publications ult\'{e}rieures d\'{e}velopperont les consid\'{e}rations
ci-dessous. Elles ont pour ambition de provoquer une discussion \og
\'{e}pist\'{e}mologique \fg, aussi ouverte que possible, en math\'{e}matiques
appliqu\'{e}es.
\subsection{Automatique}
Quel sens gardent les th\`{e}mes traditionnels de la recherche en
automatique si les th\`{e}ses de la commande sans-mod\`{e}le sont \og
vraies\footnote{On entend par \og v\'{e}rit\'{e} \fg \ la confirmation des
succ\`{e}s pratiques.}\fg? Que l'on songe par exemple, d'un point de vue
strictement acad\'{e}mique, \`{a}
\begin{itemize}
\item la plupart des questions structurelles en lin\'{e}aire
et en non-lin\'{e}aire de dimension finie, puisqu'il n'y a plus de
mod\`{e}les;
\item l'identification param\'{e}trique, car sans mod\`{e}les point de param\`{e}tres!;
\item la commande robuste et la commande stochastique,
puisque l'iPI prend en compte les termes inconnus et les
perturbations;
\item la commande optimale, d\'{e}terministe ou stochastique: le
caract\`{e}re toujours arbitraire et artificiel du choix du
crit\`{e}re\footnote{En ing\'{e}nierie, contrairement \`{a} la physique \og pure
\fg, o\`{u} les lois fondamentales d\'{e}rivent de principes variationnels,
les crit\`{e}res utilis\'{e}s pour l'optimisation sont \og bricol\'{e}s \fg.}
est accentu\'{e} par la m\'{e}connaissance du mod\`{e}le.
\end{itemize}
Aujourd'hui, l'exploration des limites du sans-mod\`{e}le semble une
voie \`{a} privil\'{e}gier, en particulier par rapport \`{a} la dimension
infinie, c'est-\`{a}-dire par rapport aux retards et aux \'{e}quations aux
d\'{e}riv\'{e}es partielles. Les succ\`{e}s obtenus, notamment avec des syst\`{e}mes
o\`{u} la mod\'{e}lisation par \'{e}quations aux d\'{e}riv\'{e}es partielles est
primordiale, comme en \cite{edf}, avec les installations
hydro\'{e}lectriques, et en \cite{agadir}, avec la r\'{e}gulation du trafic
autoroutier, rendent cette investigation urgente.

\subsection{Mod\'{e}lisation math\'{e}matique}
Une large part de la recherche universitaire appliqu\'{e}e souligne,
par-del\`{a} l'automatique, la n\'{e}cessit\'{e} de l'\'{e}criture de mod\`{e}les
math\'{e}matiques fiables, dont les \'{e}quations diff\'{e}rentielles, ou aux
diff\'{e}rences, sont l'armature. C'est trop souvent un leurre ! Y
renoncer
\begin{itemize}
\item n'implique pas l'abandon des math\'{e}matiques;
\item signifie, comme indiqu\'{e} ici, et en \cite{finance} pour la finance quantitative,
la poursuite de concepts et, donc, d'outils th\'{e}oriques nouveaux.
\end{itemize}

\end{document}